\documentclass[a4paper,12pt,makeidx]{amsart}
\usepackage{amscd}
\usepackage{xypic}  
\usepackage{amssymb}
\usepackage{amsthm}
\usepackage{epsf}
\makeindex
\newtheoremstyle{break}
  {9pt}
  {9pt}
  {\itshape}
  {}
  {\bfseries}
  {.}
  {\newline}
  {}
\theoremstyle{break}
\newtheorem{bthm}{Theorem}

\newtheorem{bcor}{Corollary}
\theoremstyle{plain}
\newtheorem{thm}{Theorem}[section]
\newtheorem{cor}[thm]{Corollary}
\newtheorem{exa}[thm]{Example}
\newtheorem{lemma}[thm]{Lemma}

\newtheorem{prop}[thm]{Proposition}

\newtheorem{defn}[thm]{Definition} 
\newtheorem{notn}[thm]{Notation} 
\newtheorem{claim}[thm]{Claim}
\newtheorem{rem}[thm]{Remark}

\renewcommand{\proofname}{Proof}

\def\Sing{\operatorname{Sing}}
\def\Ker{\operatorname{Ker}}
\def\Coker{\operatorname{Coker}}

\def\min{\operatorname{min}}
\def\Im{\operatorname{Im}}

\def\max{\operatorname{max}}
\def\length{\operatorname{length}}
\def\c1{\operatorname{c_1}}
\def\c2{\operatorname{c_2}}
\def\Cliff{\operatorname{Cliff}}

\def\Sym{\operatorname{Sym}}
\def\RR{{\mathbb R}}

\def\ZZ{{\mathbb Z}}
\def\NN{{\mathbb N}}
\def\QQ{{\mathbb Q}}
\def\PP{{\mathbb P}}
\def\GG{{\mathbb G}}

\def\A{{\mathcal A}}
\def\U{{\mathcal U}}
\def\V{{\mathcal V}}
\def\W{{\mathcal P}}

\def\L{{\mathcal L}}

\def\O{{\mathcal O}}
\def\I{{\mathcal J}}
 
\def\E{{\mathcal E}}

\def\H{{\mathcal H}}
\def\F{{\mathcal F}}

\def\M{{\mathcal M}}

\def\eqv{\equiv}

\def\+{\oplus}                   
\def\*{\otimes}                  
\def\hpil{\longrightarrow}       
\def\khpil{\rightarrow}

\def\Pic{\operatorname{Pic}}
\def\Gal{\operatorname{Gal}}

\def\Spec{\operatorname{Spec}}

\def\det{\operatorname{det}}
\def\Bs{\operatorname{Bs}}

\begin{document}  

\title[On the extendability of elliptic surfaces of rank two and higher]{On the 
extendability of elliptic surfaces of rank two and higher}

\dedicatory{\normalsize \dag \ Dedicated to the memory of Giulia Semproni} 

\author[A.F. Lopez, R. Mu\~{n}oz and J.C. Sierra]{Angelo Felice Lopez*, Roberto Mu\~{n}oz**
and Jos\'e Carlos Sierra**}

\address{\hskip -.43cm Angelo Felice Lopez, Dipartimento di Matematica, 
Universit\`a di Roma Tre, Largo San Leonardo Murialdo 1, 00146, Roma, Italy. e-mail {\tt
lopez@mat.uniroma3.it}}

\address{\hskip -.43cm Roberto Mu\~{n}oz, Departamento de Matem\'atica Aplicada, Universidad Rey Juan Carlos, 28933 M\'ostoles 
(Madrid), Spain. e-mail {\tt roberto.munoz@urjc.es}} 

\address{\hskip -.43cm Jos\'e Carlos Sierra, Departamento de \'Algebra, Facultad de Ciencias Matem\'aticas, Universidad
Complutense de Madrid, 28040 Madrid, Spain. e-mail {\tt jcsierra@mat.ucm.es}}

\thanks{* Research partially supported by the MIUR national project ``Geometria delle variet\`a algebriche"
COFIN 2002-2004 and by the INdAM project ``Geometria birazionale delle variet\`a algebriche". The author would
like to thank the CIMAT (Guanajuato, Mexico) for the wonderful hospitality during part of this research.}

\thanks{** Research partially supported by MCYT project BFM2003-03971. The third author was also partially supported by the 
program ``Profesores de la UCM en el extranjero. Convocatoria 2006".}


\begin{abstract}

We study threefolds $X \subset \PP^r$ having as hyperplane section a smooth surface with an elliptic fibration. We first
give a general theorem about the possible embeddings of such surfaces with Picard number two. More precise results are then
proved for Weierstrass fibrations, both of rank two and higher. In particular we prove that a Weierstrass fibration of rank
two that is not a K3 surface is not hyperplane section of a locally complete intersection threefold and we give some
conditions, for many embeddings of Weierstrass fibrations of any rank, under which every such threefold must be a cone.
\end{abstract}

\maketitle
 
\section{Introduction}
\label{intro}

The minimal model program has highlighted the importance, among the basic building blocks in the study of
birational equivalence classes of algebraic varieties, of Mori fiber spaces, that is morphisms $f : X
\to Y$ with connected fibers such that $X$ is normal projective with $\QQ$-factorial terminal singularities, $Y$ is normal
projective, $\dim Y < \dim X$, $-K_X$ $f$-ample and $\rho(X) - \rho(Y) = 1$. In dimension $3$, where the minimal model
program has been accomplished, a lot of work has been dedicated to the study of the case when $Y$ is a point,
that is when $X$ is a Fano variety (\cite{isk1, isk2, mm, pr}). Perhaps the next interesting case is when $Y$ is a
curve and here also several papers have appeared (see \cite{bcz} and references therein).

In the present article we also study the latter case, but from a different point of view, that we wish to outline here.
Given a Mori fiber space $f : X \to Y$ with $Y$ a curve and general fiber $F$, in many cases we can take a projective
embedding $X \subset \PP^r$ with $\O_X(1) \cong \O_X(- K_X + h F)$, $h >> 0$. Now a general hyperplane section
$S = X \cap H$ inherits an elliptic fibration and will often have $\rho(S) = \rho(X) = 2$ (the latter happens, for
example when $X$ is smooth, $h^2(\O_X) = 0$ and $p_g(S) > 0$, by a theorem of Moishezon \cite[Thm.7.5]{moi}). 

Reversing this scenery it seems therefore interesting to take a projective embedding $S \subset \PP^N$ of an elliptic surface
$S$, for example with Picard number two, and study which threefolds $X \subset \PP^{N+1}$ can have $S$ as hyperplane section.
In the literature there is also a lot of work performed in this direction, mostly based on the following two techniques. The
first one is adjunction theory (see \cite{bs}), that, to say in a few words in the specific cases we are describing, aims first at
extending the fibration to the threefold and then to study the properties of the threefold using the extended fibration. A nice
example of this is the result of Badescu (\cite[Thm.7]{ba}), that classifies pairs $(X, L)$ with $X$ a normal projective variety, $L$
an ample line bundle on $X$ such that there exists $S \in |L|$ that is a $\PP^s$-bundle over a curve. The second technique is more
recent and is based on a theorem of Zak \cite[page 277]{z} and on the theory of Gaussian maps \cite{w}. For this point of view we
mention here the study of smooth Fano threefolds and Mukai varieties \cite{clm1, clm2}, and, more recently, of Enriques-Fano
threefolds and threefolds with hyperplane sections pluricanonical surfaces of general type \cite{klm}.

To explain the results proved in this article, employing both techniques above, we start with a few definitions. 

In the sequel all varieties are over the complex numbers.

\begin{defn}
A subvariety $Y \subset \PP^N$ is called {\bf extendable} if there exists a subvariety $X \subset \PP^{N+1}$ and a
hyperplane $H = \PP^N \subset \PP^{N+1}$ such that $Y = X \cap H$, $X$ is different from a cone over $Y$ and $\dim X
= \dim Y+1$. Such an $X$ is called an {\bf extension} of $Y$.

If $Y$ is extendable to a variety $X$ as above with locally complete intersection singularities, we will
say that $Y$ is {\bf l.c.i.\ extendable}. Similarly we can define {\bf smoothly extendable} if
$X$ is smooth, or {\bf normally extendable} if $X$ is normal, {\bf l.c.i.-terminal extendable} if $X$ has terminal
locally complete intersection singularities, {\bf l.c.i.r.s.\ extendable} if $X$ is locally complete
intersection with rational singularities.
\end{defn}

We will study extensions of elliptic surfaces, as in the ensuing

\begin{defn} 
Let $S$ be a smooth irreducible projective surface and let $B$ be a smooth irreducible curve. We will denote by 
$d(B)$ the minimum degree of a very ample line bundle on $B$. An {\bf elliptic fibration} $\pi : S \to B$ is a
surjective morphism whose general fiber is a smooth connected curve of genus one. If a smooth surface $S$ has an elliptic
fibration we will call $S$ an {\bf elliptic surface}. 
\end{defn}

A simple but important point for us is that, in many cases (see Proposition \ref{mfs}), extensions of
elliptic surfaces are in fact {\bf Mori fiber spaces}.

\vskip .3cm

Our first result, which, together with Theorem \ref{fibra} below, can be considered a more precise version of
\cite[3.5.2]{mor}, studies which embeddings can occur for some extendable elliptic surfaces $S \subset \PP^N$ with
Picard number two.

\begin{bthm} 
\label{nonextell}
Let $S \subset \PP^N$ be a smooth surface having an elliptic fibration $\pi: S \to B$ with general fiber $f$. Suppose 
that $N^1(S) \cong \ZZ [C] \oplus \ZZ [f]$ for some divisor $C$ and that the hyperplane bundle of $S$ is $H_S \eqv aC +
bf$, so that we can also suppose, without loss of generality, that $C.f \geq 1$. Let $X \subset \PP^{N+1}$ be any l.c.i.\
extension of $S$ and suppose furthermore that one of the  following holds:
\begin{itemize}
\item[(i)] $X$ has rational singularities and $g(B) > 0$;
\item[(ii)] $X$ has $\QQ$-factorial terminal singularities, $B \cong \PP^1$ and $\kappa(S) = 1$;
\item[(iii)] $\kappa(S) = 1$ and $H^1(S, K_S + H_S - f_1 - \ldots - f_{d(B)}) = 0$ for every set of smooth distinct fibers
$f_i$'s.
\end{itemize}

Then 
\[ (a, C.f) \in \{ (1, 3), (1, 4), (1, 5), (1, 6), (1, 7), (1, 8), (1, 9), (2, 4), (3, 3) \}. \]
Moreover if, in addition, $X$ is locally factorial, then $(a, C.f) \not\in \{(1, 7), (1, 8), (1, 9) \}.$ 
\end{bthm}

The vanishing condition in (iii) above is satisfied in many cases, see for example Remark \ref{van}. Moreover
we observe that there are examples of smoothly extendable smooth elliptic surfaces with $\kappa(S) = 1$, $N^1(S)
\cong \ZZ [C] \oplus \ZZ [f]$, $H_S \eqv aC + bf$ and $(a, C.f) \in \{ (1, 3), (1, 4), (1, 5), (1, 6), (2, 4), (3, 3)
\}$ (see examples \ref{uno}-\ref{sei}).

\vskip .3cm

Our results can be made a lot more precise if we assume a little bit more on the fibration: In case $S$ has a Weierstrass
fibration (see Definition \ref{weifibr}) and $\rho(S) = 2$, then $S$ is often non l.c.i.\ extendable (but it can be extendable in
higher rank, see Example \ref{sette}).

\begin{bcor} 
\label{nonextweier}
Let $S \subset \PP^N$ be a smooth surface having a Weierstrass fibration $\pi: S \to B$ with general fiber $f$
and section $C$. Set $n = - C^2, g = g(B)$ and suppose that $\rho(S) = 2$, $n \geq 1$ and $(g, n) \neq (0, 1)$.
\begin{itemize}
\item[(i)] If $(g, n) \neq (0, 2)$ then $S$ is not l.c.i.\ extendable.
\item[(ii)] If $(g, n) = (0, 2)$ then any possible l.c.i.\ extension $X \subset \PP^{N+1}$ of
$S$ is an anticanonically embedded Fano threefold with $\rho(X) = 1$ and $h^1(\O_X) = h^2(\O_X) = 0$.
\end{itemize}
\end{bcor}

In the case (ii), which turns out to be exactly the $K3$-Weierstrass case (see Proposition \ref{invariants}(iii)), we
can be a little bit more precise.

\begin{bcor} 
\label{nonextweierk3}
Let $S \subset \PP^N$ be a smooth surface having a Weierstrass fibration $\pi: S \to \PP^1$ with general fiber
$f$ and section $C$ such that $C^2 = -2$. Suppose that $\rho(S) = 2$. Let $H_S \sim aC + bf$ be the hyperplane
bundle of $S$ and let $g(S)$ be the sectional genus of $S$. We have:

{\rm (i)} $S$ is not l.c.i.-terminal extendable. 

{\rm (ii)} If $(a, b, g(S)) \not\in \{(3, 7, 13), (3, 8, 16), (3, 10, 22), (3, 11, 25), (3, 13, 31), (3, 14, 34), (4, 9,
21),$ 

\ \null \hskip.5cm $(4, 11, 29), (4, 13, 37), (5, 11, 31), (5, 12, 36) \}$, then $S$ is not l.c.i.\ extendable.

{\rm (iii)} If $(a, b, g(S)) \not\in \{(3, 7, 13), (3, 8, 16), (3, 9, 19), (3, 10, 22), (3, 11, 25), (3,
12, 28), (3, 13, 31),$ 

\ \null \hskip.5cm $(3, 14, 34), (3, 15, 37), (4, 9, 21), (4, 10, 25), (4, 11, 29), (4, 12, 33), (4, 13, 37), (5,
11, 31),$ 

\ \null \hskip.48cm $(5, 12, 36) \}$, then $S$ is not normally extendable.
\end{bcor}

We end this introduction with a non extendability result (regardless of the singularities of the extension) for
Weierstrass fibrations with more special assumptions on the embedding line bundle, {\it but with no assumption on
the rank on the Picard group}.

\begin{bthm} 
\label{nonextweier2}
Let $S \subset \PP^N$ be a smooth surface having a Weierstrass fibration $\pi: S \to B$ with general fiber $f$ and section
$C$.  Set $n = - C^2$ and $g = g(B)$. Suppose that the hyperplane bundle of $S$ is of type $H_S \eqv aC + bf$ and that $n
\geq 1$. 

Then $S$ is not extendable if any of the following conditions is satisfied:
\begin{itemize}
\item[(i)] $g = 0$ and $a \geq 6$ with $(a, b, n) \neq (6, 7, 1)$, or
\item[(ii)] $g \geq 1$ and $a \eqv 0 \ ({\rm mod} \ 3)$, $a \geq 6$, or 
\item[(iii)] $g \geq 1$, $S$ is linearly normal and either $a \geq 7$, $b \geq an + 5g - 1$ or 
$a = 5$, $b \geq 6n + 7g - 3$.
\end{itemize}
\end{bthm}

\vskip .2cm

\noindent {\it Acknowledgments}. The authors wish to thank Mike Roth for several helpful discussions.

\section{Background material}
\label{backgr}

We recall in this section, for the reader's convenience, some results of adjunction theory that will be used in the sequel.

Throughout this section we will denote by $X$ an irreducible normal $n$-dimensional projective variety with terminal singularities, $n
\geq 2$, by $A$ an ample line bundle on $X$ and by $r$ the index of singularities of $X$, that is the smallest positive integer $r$
such that $rK_X$ is a Cartier divisor.

\begin{defn} \cite[Def.1.5.3]{bs} 
\label{nefval}
The {\bf nefvalue of $(X, A)$} is 
\[\tau (X, A) = \min \{t \in \RR : K_X + tA \ \mbox{is nef} \}. \]
\end{defn}

By Kawamata's rationality theorem \cite[Thm.4.1.1]{kmm}, if $K_X$ is not nef, the nefvalue is a rational number and we can write $r
\tau = \frac{u}{v}$ for some $u, v \in \NN$. By the Kawamata-Shokurov base-point free theorem \cite[Thm.3.3]{km}, the linear system
$|m(vrK_X + uA)|$ is  base-point free for $m >> 0$ and using the Stein factorization of the morphism defined by this linear system one
gets a morphism $\phi(X, A) : X \to Y$ with connected fibers onto a normal projective variety $Y$. This morphism depends only on the
pair $(X, A)$ and is called the {\bf nefvalue morphism of $(X, A)$}.

We mention here several useful results of general adjunction theory related to the nefvalue.

\begin{prop} \cite[Prop.7.2.2]{bs} 
\label{7.2.2}
Let $\tau$ be nefvalue of $(X, A)$. Then either
\begin{itemize}
\item[(i)] $\tau = n + 1$ and $(X, A) \cong (\PP^n, \O_{\PP^n}(1))$; or
\item[(ii)] $\tau = n$ and $(X, A) \cong (Q, \O_Q(1))$, $Q \subset \PP^{n+1}$ a quadric; or
\item[(iii)] $\tau = n$ and $(X, A)$ is a $(\PP^{n-1}, \O_{\PP^{n-1}}(1))$-bundle over a smooth curve under $\phi(X, A)$; or
\item[(iv)] $\tau \leq n$ and $K_X + nA$ is nef and big.
\end{itemize}
\end{prop}

In the range $n - 1 < \tau < n$ we have

\begin{thm} \cite[Thms.7.2.3 and 7.2.4]{bs} 
\label{7.2.34}
Assume that $X$ is $\QQ$-factorial and suppose that $K_X + nA$ is nef and big. Then $K_X + nA$ is ample and if $\tau$ is
the nefvalue of $(X, A)$ we have $\tau \leq n - 1$ unless $\tau = n - \frac{1}{2}$ and $(X, A)$ is a generalized cone over $(\PP^2,
\O_{\PP^2}(2))$.
\end{thm} 

When $\tau \leq n - 1$ we have

\begin{thm} \cite[Thm.7.3.2]{bs} 
\label{7.3.2}
Suppose that $\tau (X, A) \leq n - 1$. Then $K_X + (n - 1)A$ is ample unless $\tau (X, A) = n - 1$ and either
\begin{itemize}
\item[(i)] $rK_X \sim -r(n - 1)A$; or
\item[(ii)] $(X, A)$ is a quadric fibration over a smooth curve under $\phi(X, A)$; or
\item[(iii)] $(X, A)$ is a scroll over a normal surface under $\phi(X, A)$; or
\item[(iv)] $\phi = \phi(X, A) : X \to Y$ is birational. Moreover if $X$ is factorial then $\phi(X, A)$ is the simultaneous
contraction to distinct smooth points of divisors $E_i \cong \PP^{n-1}$ such that $E_i \subset {\rm Reg}(X)$, $\O_{E_i}(E_i) \cong
\O_{\PP^{n-1}}(-1)$ and $A_{|E_i} \cong \O_{\PP^{n-1}}(1)$.  Also $L := (\phi_{\ast}(A))^{\ast \ast}$ and $K_Y + (n - 1)L$ are 
ample and $K_X + (n - 1)A \cong \phi^{\ast}(K_Y + (n - 1)L)$.
\end{itemize}
\end{thm}

When $K_X + (n - 1)A$ is nef and big we can define the first reduction of $(X, A)$. By the Kawamata-Shokurov base-point free theorem
\cite[Thm.3.3]{km} we have a birational morphism $\pi : X \to X'$ with connected fibers and normal image associated to the linear
system $|mr(K_X + (n-1)A)|$ for $m >> 0$. Set $A' := (\pi_{\ast}(A))^{\ast \ast}$. The pair $(X', A')$ is called the {\bf
first reduction of $(X, A)$}. Then we have

\begin{thm} \cite[Thm.7.3.4]{bs} 
\label{7.3.4}
Assume that $X$ is factorial and that $n \geq 3$ and let $(X', A')$ be the first reduction of $(X, A)$. Suppose that $n - 2 < \tau
(X', A') < n - 1$. Then either
\begin{itemize}
\item[(i)] $n = 4, \tau (X', A') = \frac{5}{2}$ and $(X', A') \cong (\PP^4, \O_{\PP^4}(2))$; or
\item[(ii)] $n = 3, \tau (X', A') = \frac{3}{2}$ and $(X', A') \cong (Q, \O_Q(2))$, $Q \subset \PP^4$ a quadric; or
\item[(iii)] $n = 3, \tau (X', A') = \frac{4}{3}$ and $(X', A') \cong (\PP^3, \O_{\PP^3}(3))$; or
\item[(iv)] $n = 3, \tau (X', A') = \frac{3}{2}$, $\phi(X, A)$ has a smooth curve as image and $(F, A'_{|F}) \cong (\PP^2,
\O_{\PP^2}(2))$ for a general fiber $F$ of $\phi(X, A)$.
\end{itemize}
\end{thm}

\section{Weierstrass fibrations}

We collect in this section some notation and facts about Weierstrass fibrations that will be used in the sequel.

Let $\pi: S \to B$ be a minimal elliptic surface, that is there are no $(-1)$-curves in the fibers of
$\pi$. Suppose moreover that $\pi$ has a section $s : B \to S$. On each reducible fiber contract any irreducible
component not meeting $s(B)$. Hence we obtain a new (singular) elliptic surface $\pi': S' \to B$ with a
section and whose fibers are all both reduced and irreducible. In this context a global Weierstrass equation can
be given and the following concept appears (see \cite{mi2}).

\begin{defn}
\label{weifibr}
Let $S$ be a surface and let $B$ be a smooth curve. A {\bf Weierstrass fibration} $\pi: S \to B$ is a flat and
proper map such that every geometric fiber has arithmetic genus one (so that it is either a smooth genus one curve,
or a rational curve with a node, or a rational curve with a cusp), with general fiber smooth and such that there
is given a section of $\pi$ not passing through the singular point of any fiber.

We will say that a Weierstrass fibration $\pi: S \to B$ is {\bf smooth} if $S$ is smooth. 
\end{defn}

\begin{rem}
{\rm By the above discussion the notions of {\bf Weierstrass fibration} and {\bf elliptic surface with section} can
be freely interchanged when $S$ is smooth and $\rho(S) = 2$.}
\end{rem}

We recall from \cite{mi2} some well-known facts about Weierstrass fibrations (see also \cite{mi1}).

\begin{defn}
Let $\pi: S \to B$ be a Weierstrass fibration with section $C \subset S$. We define the {\bf fundamental line
bundle} $\L$ of $\pi$ as the dual line bundle of $\pi_{\ast} N_{C/S}$ on $B$. We will set $n = \deg \L$.
\end{defn}

\begin{rem} 
{\rm As $\L = (R^1\pi_{\ast} \O_S)^{\ast}$, the fundamental line bundle $\L$ does not depend on the given section
$C$. Moreover $n = -C^2 \geq 0$ and $\L \cong \O_B$ if and only if $S$ is a product $B \times F$, with $F$ an
elliptic curve \cite[(II.3.6) and (III.1.4)]{mi2}.}
\end{rem}

\begin{lemma} \cite[(II.3.5), (II.3.7) and (II.4.3)]{mi2}
\label{split}
Let $\pi: S \to B$ be a Weierstrass fibration with section $C$ and fundamental line bundle $\L$. We have:
\begin{itemize}
\item[(i)] $\pi_{\ast} \O_S \cong \pi_{\ast} \O_S(C) \cong \O_B$, $R^1 \pi_{\ast} \O_S(uC) = 0$ for every $u
\geq 1$.
\item[(ii)] $\pi_{\ast} \O_S(mC) \cong \O_B \oplus \L^{-2} \oplus \L^{-3} \oplus \ldots \oplus\L^{-m} \mbox{ for
every } m \geq 2.$
\end{itemize}
\end{lemma}

The invariants and Kodaira type of a smooth Weierstrass fibration are given as follows.

\begin{prop} \cite[(III.1.1), (III.4), (IV.1.1) and (VII.1.3)]{mi2}
\label{invariants}
Let $\pi: S \to B$ be a smooth Weierstrass fibration with section $C$, general fiber $f$ and fundamental line
bundle $\L$ with $n = \deg \L \geq 1$ and $g = g(B)$. We have: 
\begin{itemize}
\item[(i)] $K_S \eqv (n + 2g - 2) f$, $q(S) = h^1(\O_S) = g$, $p_g(S) = h^0(K_S) = n - 1 + g$ and
$h^{1,1}(S) = 10 n + 2g$.
\item[(ii)] $\kappa(S) = -\infty$ if and only if $S$ is a rational surface if and only if $g = 0$ and
$n = 1$.
\item[(iii)] $\kappa(S) = 0$ if and only if $S$ is a $K3$ surface if and only if $g = 0$ and
$n = 2$. 
\item[(iv)] $\kappa(S) = 1$ if and only if $(g, n) \not\in \{(0, 1), (0, 2) \}$.
\item[(v)] Let $A \in \Pic S$. Then $A \eqv \alpha C + \beta f$ if and only if $A \cong \O_S(\alpha C) \otimes
\pi^{\ast} M$, for some $M \in \Pic^{\beta} B$.
\end{itemize}
\end{prop}

By the above Proposition, the rank of the Picard group of a smooth Weierstrass fibration satisfies
\[ 2 \leq \rho(S) \leq h^{1,1}(S) = 10 n + 2g. \]
On the other hand, for surfaces with $p_g > 0$ the Picard number $\rho$ is in general strictly less than
$h^{1,1}$ and by Hodge theory one expects at least $p_g$ independent conditions for a given cycle to be algebraic.
The typical picture one conjectures is that the generic surface in its moduli space has low Picard number (as it
is for general surfaces in $\PP^3$).

For Weierstrass fibrations over $\PP^1$ this prediction turns out to be true. In the latter case Miranda \cite{mi1}
constructed a moduli space for such fibrations and Cox \cite[MainThm.]{co} (see also \cite[Cor.1.2]{klo}) proved
that a general (in the countable Zariski topology) Weierstrass fibration $\pi : S \to \PP^1$ with $n \geq 2$,
satisfies $\rho(S) = 2$.

\begin{rem}
\label{weirk2}
It is likely that there exist smooth Weierstrass fibrations $\pi: S \to B$ with $\rho(S) = 2$, where 
$B$ is a smooth curve of genus $g$, for every $g \geq 1$. According to a suggestion of R. Kloosterman they
should be constructed as follows. Let $E$ be the elliptic curve with J-invariant zero and, for any $n \geq 1$, let
$P_1, \ldots, P_{6n} \in B$ and let $C$ be a cyclic covering of degree $6$ ramified at $P_1, \ldots, P_{6n}$.
If $C$ does not have a nonconstant morphism to $E$ then $\rho(S) = 2$. {\rm To see this note that there is a
line bundle $\L$ of degree $n$ on $B$ such that $\L^6 \cong \O_B(P_1 + \ldots + P_{6n})$. We have
therefore a nonzero section $s \in H^0(\L^6)$ giving rise to the Weierstrass data $(\L, 0, s)$ on $B$ and, by
\cite[II.5]{mi2}, to a smooth Weierstrass fibration $\pi: S \to B$. By ``Tate's algorithm"
\cite[IV.3.1]{mi2} all fibers of $\pi$ are irreducible (the singular ones being cuspidal). By Shioda-Tate's formula
\cite[Cor.VII.2.4]{mi2}, we have that $\rho(S) = 2$ if the Mordell-Weil group of sections of $\pi$ has rank zero, that is if all
sections are of finite order. But a section of infinite order gives, as in \cite[section 6]{klo}, a nonconstant morphism from $C$
to the elliptic curve $E$ with J-invariant zero.}
\end{rem}

The following remark will be useful.

\begin{rem}
\label{rho=2} 
Let $\pi: S \to B$ be a smooth Weierstrass fibration with section $C$ and general fiber $f$. If $\rho(S) = 2$ then $N^1(S) \cong \ZZ
[C] \oplus \ZZ [f]$. {\rm To see this note that, since $N^1(S)$ is torsion free, we have $N^1(S) \cong \ZZ [A] \oplus \ZZ [A']$ for
some $A, A' \in \Pic S$. On the other hand $\ZZ [C] \oplus \ZZ [f]$ has rank two, therefore there are integers $a \geq 1, a' \geq 1,
u, v, u', v'$ such that $aA \eqv uC + vf$ and $a'A' \eqv u'C + v'f$. Now $a A.f = u$, whence $a$ divides both $u$ and $v$ and we
get $A \eqv u_1C + v_1f$ for some integers $u_1, v_1$. Similarly $A' \eqv u_1'C + v_1'f$.}
\end{rem}

To study the extendability of Weierstrass fibrations we will need the following simple results.

\begin{lemma}
\label{van1} 
Let $\pi: S \to B$ be a smooth Weierstrass fibration with section $C$, general fiber $f$ and
fundamental line bundle $\L$ with $n = \deg \L \geq 1$ and $g = g(B)$. Let $D \eqv \alpha C + \beta f$. For $g
\geq 1$ and $\W \in \Pic^0 B$ we will set $D_{\W} = D \otimes \pi^{\ast} \W$. We have:
\begin{itemize}
\item[(i)] $H^1(D) = 0$ if either $\alpha = 1$ and $\beta \geq 2g - 1$ or $\alpha \geq 2$ and $\beta  \geq \alpha 
n + 2g - 1$.
\item[(ii)] If $g \geq 1$ and $\W \in \Pic^0 B$ is general, then $H^1(D_{\W}) = 0$ if either $\alpha = 1$ and
$\beta \geq g - 1$ or $\alpha \geq 2$ and $\beta  \geq \alpha n + g - 1$.
\item[(iii)] If $H \eqv aC + b f$, $g \geq 1$ and $\W \in \Pic^0 B$ is general, then $H^1(H - 2D_{\W}) = 0$ if
either $a - 2\alpha = 1$ and $b - 2\beta \geq g - 1$ or $a - 2\alpha \geq 2$ and $b - 2\beta  \geq \alpha n + g -
1$.
\item[(iv)] $|D|$ is base-point free if $\alpha \geq 2$ and $\beta \geq \alpha n + 2g$.
\item[(v)] $D$ is very ample if $\alpha \geq 3$ and $\beta \geq \alpha n + 2g + 1$.
\item[(vi)] If $D$ is very ample then $\alpha \geq 3$ and $\beta \geq \alpha n + 1$.
\end{itemize}
\end{lemma}

\begin{proof} To see (i) note that, by Proposition \ref{invariants}(v), we have $D \cong \O_S(\alpha C)
\otimes \pi^{\ast} M$ for some $M \in \Pic^{\beta} B$. By Lemma \ref{split}(i) and the Leray spectral sequence we
deduce that $H^1(S, D) \cong H^1(B, \pi_{\ast} D) = 0$ for degree reasons.

We now show (ii) and (iii). Let $d \geq g - 1$ be an integer. Since $g \geq 1$ we know that there is a nonempty
open subset $\V_d \subset \Pic^d B$ such that $H^1(L) = 0$ for any $L \in \V_d$. Given any $N \in \Pic^d B$
consider the isomorphism $\phi_N : \Pic^0 B \to \Pic^d B$ given by tensoring with $N$. Then we get a nonempty open
subset $\U_N := \phi_N^{-1}(\V_d) \subset \Pic^0 B$ such that $H^1(N \otimes \W) = 0$ for any $\W \in \U_N$.

As above we have $D \cong \O_S(\alpha C) \otimes \pi^{\ast} M$ for some $M \in \Pic^{\beta} B$ and $D_{\W} \cong
\O_S(\alpha C) \otimes \pi^{\ast} (M \otimes \W)$. Since $\alpha \geq 1$, Lemma \ref{split}(i) and the Leray
spectral sequence imply that $H^1(S, D_{\W}) \cong H^1(B, \pi_{\ast} D_{\W})$. By Lemma \ref{split}(i) and (ii) we
have $\pi_{\ast} D_{\W} = M \otimes \W$ when $\alpha = 1$ and $\pi_{\ast} D_{\W} = (M \otimes \W) \oplus
\bigoplus\limits_{i=2}^{\alpha} (M \otimes \L^{-i} \otimes \W)$ when $\alpha \geq 2$. Hence, to prove (ii), we can
choose $\W \in \U_D := \U_M$ when $\alpha = 1$ and $\W \in \U_D := \U_M \cap \bigcap\limits_{i=2}^{\alpha} \U_{M
\otimes \L^{-i}}$ when $\alpha \geq 2$.

Now to see (iii) consider, using additive notation for line bundles on $B$, the surjective morphism $h_2 : \Pic^0
B \to \Pic^0 B$ defined by $h_2(\W) = -2 \W$. Hence, given any $N \in \Pic^d B$, we get a surjective morphism 
$\psi_N := \phi_N \circ h_2 : \Pic^0 B \to \Pic^d B$. Therefore $H^1(N - 2 \W) = 0$ for any $\W \in \A_N := \psi_N^{-1}
(\V_d)$. Now by Proposition \ref{invariants}(v), we have
$H \cong \O_S(a C) \otimes \pi^{\ast} M_H$ for some $M_H \in \Pic^b B$, whence $H - 2D_{\W} \cong \O_S((a - 2
\alpha) C) \otimes \pi^{\ast} (M_H - 2M - 2 \W)$. Since $a - 2 \alpha \geq 1$, Lemma \ref{split}(i) and the Leray
spectral sequence imply that $H^1(S, H - 2D_{\W}) \cong H^1(B, \pi_{\ast} (H - 2D_{\W}))$. By Lemma \ref{split}(i)
and (ii) we have $\pi_{\ast} (H - 2D_{\W}) = M_H - 2M - 2 \W$ when $a - 2 \alpha = 1$ and $\pi_{\ast} (H -
2D_{\W}) = (M_H - 2M - 2 \W) \oplus \bigoplus\limits_{i=2}^{a - 2 \alpha} (M_H - 2M - i\L - 2 \W)$ when $a - 2 \alpha
\geq 2$. Hence, to prove (iii), we can choose $\W \in \U_{H, D} := \A_{M_H - 2M}$ when $a - 2 \alpha = 1$ and $\W
\in \U_{H, D} := \A_{M_H - 2M} \cap \bigcap\limits_{i=2}^{a - 2 \alpha} \A_{M_H - 2M - i\L}$ when $a - 2 \alpha \geq 2$.

To prove (iv) note that $H^1(D - f) = 0$ by (i). Let $F$ be any fiber of $\pi$. We will be done if we prove
that $|D_{|F}|$ is base-point free. Now this follows by \cite[Prop.2.3,I]{cf} since we have $\alpha = D.F \geq 2 =
2p_a(F)$.

Similarly to see (v) note that, for any fiber $F$, we have $H^1(D - F) = 0$ by (i) and $|D - F|$ is base-point free
by (iv). Let $x, y \in S$ be two distinct points. If $x$ and $y$ belong to the same fiber $F$, we can separate them
with sections in $|D|$ since $|D_{|F}|$ is very ample by \cite[Thm.3.1]{cf} (because we have $\alpha = D.F \geq 3 =
2p_a(F) + 1$). If $x$ and $y$ belong to two different fibers $F_x$ and $F_y$ respectively, then to separate them
just use the fact that $|D - F_x|$ is base-point free. On the other hand suppose that $x \in S$, $y \in T_x S$ and $d
\varphi_D (y) = 0$, where $d \varphi_D$ is the differential of the morphism $\varphi_D : S \to \PP H^0(D)$.
Arguing as above we deduce that $y$ must be tangent to $F_x$, contradicting the fact that $|D_{|F_x}|$ is very ample.

Finally (vi) is a consequence of the fact that $\alpha = D.f$ and $\beta - \alpha n = D.C$.
\end{proof}

\begin{lemma}
\label{nonhyper}
Let $\pi: S \to B$ be a smooth Weierstrass fibration with section $C$, general fiber $f$ and fundamental line
bundle $\L$ with $n = \deg \L \geq 1$ and $g = g(B)$. Let $D_0 \in \Pic S$ such that either $D_0 \eqv 3C +
\beta f$ with $\beta \geq 3n$ if $g = 0$ or $D_0 \eqv 2C + \beta f$ with $\beta \geq 2n + 2g$ if $g \geq 1$.
Then a general curve $D \in |D_0|$ is smooth irreducible and nonhyperelliptic. Moreover, if $g \geq 1$, then $D$ is
nontrigonal.

\end{lemma}
\begin{proof}
By Lemma \ref{van1}(iv) we know that $|D_0|$ is base-point free, whence $D$ is smooth and irreducible by Bertini's
theorems.

First consider the case $g = 0$. Note that $f_{|D}$ gives a $g^1_3$ on $D$. If $D$ is hyperelliptic then $D$ has
a morphism $D \to \PP^1 \times \PP^1$ which is (necessarily) birational onto its image $\overline D$. Therefore we
get the contradiction
\[ 4 \leq 6n - 2 \leq 3 \beta - 3n - 2 = g(D) \leq p_a({\overline D}) = 2. \]
Next consider the case $g \geq 1$, so that $D_0^2 = 4 \beta - 4n \geq 12$, with equality only if $n = g = 1$ and $\beta = 4$.  

If equality holds and $D$ is trigonal, using the $2:1$ morphism $\pi_{|D} : D \to B$, we get a morphism
$D \to B \times \PP^1$ which is (necessarily) birational onto its image $\overline D$. But this gives the
contradiction $8 = g(D) \leq p_a({\overline D}) = 4$.

To deal with the remaining cases, let $A$ be a base-point free $g^1_k$ on $D$, with $k = 2, 3$. By the above, we know 
that it cannot be $k = 3, n = g = 1$ and $\beta = 4$.

Let $\F = \Ker \{H^0(A) \otimes \O_S \to A \}$ and define $\E = \F^{\ast}$. As is well known (\cite{la1}), $\E$ is a
rank two vector bundle sitting in an exact sequence
\begin{equation}
\label{laz}
0 \hpil H^0(A)^{\ast} \otimes \O_S \hpil \E \hpil N_{D/S} \otimes A^{-1} \hpil 0
\end{equation}
and moreover $c_1(\E) = D$ and $c_2(\E) = k$, so that $\Delta(\E) := c_1(\E)^2 - 4c_2(\E) = D^2 - 4k > 0$.
Therefore $\E$ is Bogomolov unstable (\cite{la1}), so that, if $M$ is the maximal destabilizing subbundle (with
respect to some fixed ample line bundle $H$ on $S$), we have an exact sequence
\begin{equation}
\label{dest}
0 \hpil M \hpil \E \hpil \I_{Z/S} \otimes L \hpil 0
\end{equation}
where $L$ is another line bundle on $S$ and $Z$ is a zero-dimensional subscheme of $S$. 

We now claim that these line bundles satisfy:
\begin{itemize}
\item[(i)] $D \sim M + L$;
\item[(ii)] $k = M.L + \length(Z) \geq M.L \geq L^2 \geq 0$;
\item[(iii)] there exists an effective divisor $Z_1$ on $C$ of degree $M.L + L^2 - k \geq 0$ such that $A \cong
L_{|D}(- Z_1)$;
\item[(iv)] $L$ is base-component free and nontrivial;
\item[(v)] if $L^2 = 0$ then $M.L = k$ and $A \cong L_{|D}$.
\end{itemize} 
To see this claim note that computing Chern classes in (\ref{dest}) we get (i) and the equality in (ii). Since the
destabilizing condition reads $(M - L).H \geq 0$ and since $(M - L)^2 = \Delta(\E) + 4 \length(Z) > 0$, we
see that $M - L$ belongs to the closure of the positive cone of $S$. 

We want to prove that $\E$ is globally generated off a finite set. 

To this end observe that we just need to prove that $h^0(N_{D/S} \otimes A^{-1}) \geq 2g + 1 = 2 h^1(\O_S) + 1$,
since then the map $\psi: H^0(\E) \to H^0(N_{D/S} \otimes A^{-1})$ in (\ref{laz}) is nonzero and this gives that
$\E$ is globally generated off a finite set. Now the exact sequence
\[ 0 \hpil \O_S \hpil \O_S(D)\hpil N_{D/S} \hpil 0 \]
shows that $h^0(N_{D/S} \otimes A^{-1}) \geq h^0(N_{D/S}) - k = g - 1 + h^0(\O_S(D)) - k$, since $H^1(\O_S(D)) = 0$
by Lemma \ref{van1}(i). Using Lemma \ref{split} we find $h^0(\O_S(D)) = 2\beta - 2n - 2g + 2$, so that the desired
inequality is satisfied and $\E$ is globally generated off a finite set.

Since $\E$ is globally generated off a finite set then so is $L$. It follows that $L \geq 0$, $L$ is 
base-component free and $L^2 \geq 0$. Now the signature theorem \cite[VIII.1]{bpv} implies that $(M - L).L \geq
0$ thus proving (ii). To see (iii) and (iv) note that if $M.L > 0$ then the nefness of $L$ implies that 
$H^0(- M) = 0$.  On the other hand if $M.L = 0$ then $L^2 = D.L = 0$, whence $L \eqv 0$ by the Hodge index theorem
and therefore $D \eqv M$. Then $M.H = D.H > 0$, whence again $H^0(-M) = 0$. Twisting (\ref{laz}) and 
(\ref{dest}) by $-M$ we deduce that $h^0(L_{|D} \otimes A^{-1}) \geq h^0(\E(-M)) \geq 1$. This proves (iii) and
(iv). Moreover it gives $\deg(L_{|D} \otimes A^{-1}) \geq 0$, whence, if $L^2 = 0$, we get that $M.L \geq k$.
By (ii) it follows that $M.L = k$ and therefore $\deg(L_{|D} \otimes A^{-1}) = 0$, whence $L_{|D} \cong A$ and also
(v) is proved. 

By the Hodge index theorem we now have
\begin{equation}
\label{hodge}
L^2 (D^2 - 4k) \leq L^2 (M - L)^2 \leq (L.(M - L))^2 = (M.L - L^2)^2
\end{equation}
and it is easily seen that \eqref{hodge} gives $L^2 \leq 1$ and that $L^2 = 1$ holds precisely when $k = 3, g
= 1$ and either $n = 2, \beta = 6$ or $n = 1, \beta = 5$ (recall that we have excluded the case $k = 3, n = g = 1$
and $\beta = 4$). Moreover, in both cases above with $L^2 = 1$, we have equality in \eqref{hodge}, whence $M \eqv
3L$ and $D \eqv 4L$ by (i) above. On the other hand, in both cases, $2 = f.D$ is not divisible by $4$.

Therefore $L^2 = 0, M.L = k$ and $A \cong L_{|D}$ by (v) above. Hence $L.D = k$ and, since $L$ is nef (by (iv)
above), we must have $L.f = 0$, whence $L.C = 1$ and $k = 2$. But (iv) above also gives that $L$ is effective,
whence $L \eqv f$. From the exact sequence
\[ 0 \hpil L - D \hpil L \hpil L_{|D} \hpil 0 \]
and Lemma \ref{van1}(i) we get that $h^0(D, L_{|D}) = h^0(S, L)$. By Proposition \ref{invariants}(v) and Lemma
\ref{split}(i) we also know that $h^0(S, L) = h^0(B, \M)$ for some line bundle $\M$ of degree $1$ on $B$. Since $g
\geq 1$ we get the contradiction
\[ 2 =  h^0(A) = h^0(D, L_{|D}) = h^0(S, L) = h^0(B, \M) \leq 1. \] 
\end{proof}

\section{Extending the morphism to the threefold}
\label{extmor}

An elliptic surface has a surjective morphism onto a smooth curve with fibers elliptic curves. The goal of this
section will be to give some sufficient conditions, both on $S$ and on the singularities of $X$, to insure that
this morphism extends to a threefold containing the elliptic surface as an ample divisor. 

Extendability results for morphisms abound in the literature. We reproduce here the one of \cite[Thm.5.2.1]{bs} in
a form that will be convenient for us.

\begin{prop} 
\label{estensgener} 
Let $X$ be a projective irreducible threefold with Cohen-Macaulay singularities, let $L$ be an ample
line bundle on $X$ and let $A \in |L|$ be a normal divisor. Suppose that the restriction map $\Pic X \to \Pic A$ is an
isomorphism and that there is a surjective morphism $p: A \to Y$ onto a projective variety $Y$ such that $\dim Y \leq
\dim A - 1$.

If there exists a very ample line bundle $\L$ on $Y$ such that $H^1(A, p^{\ast} \L -
tL_{|A}) = 0$ for every $t \geq 1$, then $p$ extends to a morphism ${\overline p} : X \to Y$.
\end{prop}

\begin{proof}  
By hypothesis there is a line bundle $\H \in \Pic X$ such that $\H_{|A} \cong p^{\ast} \L$. We claim that the natural
restriction map $H^0(X, \H) \to H^0(A, \H_{|A})$ is surjective.

To this end it is of course enough to prove that $H^1(\H - L) = 0$. Now, for each $t \geq 1$, we have an exact
sequence 
\[ 0 \hpil \H - (t + 1) L \hpil \H - t L \hpil p^{\ast} \L - tL_{|A} \hpil 0, \]
whence, by hypothesis, we have that $h^1(\H - t L) \leq h^1(\H - (t + 1) L)$ for every $t \geq 1$. Let
$\omega^0_X$ be a dualizing sheaf for $X$. Since $h^1(\H - j L) = h^2(\omega^0_X \otimes (- \H + j L)) = 0$ for large
$j$ by Serre vanishing, we get that $h^1(\H - t L) = 0$ for every $t \geq 1$.

Therefore $H^0(X, \H) \to H^0(A, \H_{|A})$ is surjective, whence, since $\H_{|A}$ is globally generated, we get
$A \cap \Bs |\H| = \emptyset$. Let $m = \dim Y$ and choose, for $1 \leq i \leq m+1$, $\Delta_i \in |\L|$ such that
$\Delta_1 \cap \ldots \cap \Delta_{m+1} = \emptyset$. Pulling back to $A$ and using the above surjection we therefore
find divisors $D_i \in |\H|$ such that $D_1 \cap \ldots \cap D_{m+1} \cap A = \emptyset$. Since $A$ is ample we
have that either $D_1 \cap \ldots \cap D_{m+1} = \emptyset$ or $\dim D_1 \cap \ldots \cap D_{m+1} = 0$. But in
the latter case we get the contradiction $0 = \dim D_1 \cap \ldots \cap D_{m+1} \geq \dim X - m - 1 = \dim A - m
\geq 1$. Therefore $D_1 \cap \ldots \cap D_{m+1} = \emptyset$ and $|\H|$ is base-point free and defines a morphism
${\overline p} : X \to {\overline p}(X) \subset \PP H^0(\H)$ such that ${\overline p}_{|A} = p$. Let us show that 
${\overline p}(X) = Y$. Of course we have $Y = p(A) = {\overline p}(A) \subseteq {\overline p}(X)$. On the other hand
if there exists $x_0 \in X$ such that ${\overline p}(x_0) \not\in Y$ then $A \cap {\overline p}^{-1}({\overline
p}(x_0)) = \emptyset$, whence, as $A$ is ample, we get $\dim {\overline p}^{-1}({\overline p}(x_0)) = 0$, and
therefore $\dim {\overline p}^{-1}({\overline p}(x)) = 0$ for a general $x \in X$. Hence $\dim {\overline p}(X) =
\dim X$. Now $D_1 \cap \ldots \cap D_{m+1} = \emptyset$ implies that there are hyperplanes $H_i$ in $\PP H^0(\H)$
such that $H_1 \cap \ldots \cap H_{m+1} \cap {\overline p}(X) = \emptyset$, whence $\dim X = \dim {\overline p}(X)
\leq m \leq \dim A - 1 \leq \dim X - 2$, a contradiction.
\end{proof}

Here is an effective way to apply this to elliptic surfaces.

\begin{cor} 
\label{corestenell}
Let $X$ be a projective irreducible threefold with Cohen-Macaulay singularities, let $L$ be a very ample
line bundle on $X$ and let $S \in |L|$ be a smooth surface. Suppose that the restriction map $\Pic X \to \Pic S$
is an isomorphism and that $\pi: S \to B$ is an elliptic fibration. 

Suppose furthermore that $H^1(S, K_S + L_{|S} - f_1 - \ldots - f_{d(B)}) = 0$ for every set of distinct smooth
fibers $f_i$'s. 
  
Then $\pi$ extends to a morphism ${\overline \pi} : X \to B$. 
\end{cor}

\begin{proof}
Set $d = d(B)$ and $F_1 = 0$, $F_s = f_1 + \ldots + f_{s-1}$ for $2 \leq s \leq d + 1$. For each $t \geq 1$ and $u$
such that $1 \leq u \leq d$ we have exact sequences
\begin{equation}
\label{seq1}
0 \hpil K_S + tL_{|S} - F_{u+1} \hpil K_S + tL_{|S} - F_u \hpil (K_S + tL_{|S} - F_u)_{|f_u} \hpil 0.
\end{equation}
We first claim that $H^1(K_S + L_{|S} - F_{s+1}) = 0$ for $0 \leq s \leq d$. Since the latter is true by
hypothesis when $d = s$, we proceed by induction on $d - s$. Now when $d -s \geq 1$ we have that $H^1(K_S + L_{|S}
- F_{s+2}) = 0$ by the inductive hypothesis and $H^1((K_S + L_{|S} - F_{s+1})_{|f_{s+1}}) = 0$, whence
\eqref{seq1} with $t = 1$ and $u = s+1$ gives the claim.

Now let $W = \Im \{H^0(S, L_{|S}) \to H^0(f_s, L_{|f_s}) \}$. Then $W$ is very ample, whence $\dim W
\geq 3$. Hence we deduce, by \cite[Thm.4.e.1]{gr}, that for any $N \in \Pic B$ and for any $t \geq 1$, the
multiplication maps $\mu_{t,s} : W \otimes H^0((K_S + tL_{|S} + \pi^{\ast}N)_{|f_s}) \to H^0((K_S + (t+1)L_{|S} +
\pi^{\ast}N)_{|f_s})$ are all surjective. 

For $t \geq 1$ and $1 \leq s \leq d$ consider the restriction maps
\[ \varphi_{t,s} : H^0(K_S + tL_{|S} - F_s) \to H^0((K_S + tL_{|S} - F_s)_{|f_s}). \]
We want to prove, by induction on $t$, that they are all surjective. For $t = 1$ this follows from \eqref{seq1}
with $u = s$ and the claim proved above. Now the commutative diagram
\[ \xymatrix{
H^0(K_S + (t-1)L_{|S} - F_s)  \otimes H^0(L_{|S}) \ar[r] \ar[d]^{\varphi_{t-1,s} \* r_{f_s}}
& H^0(K_S + tL_{|S} - F_s)  \ar[d]^{\varphi_{t,s}} 
\\ H^0((K_S + (t-1)L_{|S} - F_s)_{|f_s}) \otimes W \ar[r]^{\ \ \ \ \ \ \mu_{t,s}} & H^0((K_S + tL_{|S} -
F_s)_{|f_s})    } \]
shows, by induction on $t$, that $\varphi_{t,s}$ is surjective for every $t \geq 1$. 

We now claim that $H^1(K_S + tL_{|S} - F_{s+1}) = 0$ for $t \geq 1$ and $0 \leq s \leq d$. Since the latter is
true by Kodaira vanishing when $s = 0$, we proceed by induction on $s$. Now when $s \geq 1$ we have that
$H^1(K_S + tL_{|S} - F_s) = 0$ by the inductive hypothesis, whence \eqref{seq1} with $u = s$ and the
surjectivity of $\varphi_{t,s}$ gives the claim. 

By definition of $d$ there is a very ample line bundle $\L$ of degree $d$ on $B$ and we can obviously write
$\pi^{\ast} \L \sim f_1 + \ldots + f_d$. Therefore $H^1(S, \pi^{\ast} \L - tL_{|S}) = 0$ for every $t \geq 1$ by
Serre duality and the last claim. Hence $\pi$ extends to a morphism ${\overline \pi} : X \to Y$ by Proposition
\ref{estensgener}.
\end{proof}

\begin{rem} 
\label{van}
The vanishing condition in the above corollary is satisfied in the case of smooth Weierstrass fibrations with $n \geq 1$ 
and in several other cases, even when $N^1(S) \cong \ZZ [C] \oplus \ZZ [f]$ with $C.f \geq 2$. {\rm Suppose for example 
that $\pi: S \to B$ is a smooth Weierstrass fibration with section $C$, general fiber $f$ and fundamental line bundle 
$\L$ with $n = \deg \L \geq 1$ and $g = g(B)$. Suppose furthermore that $\rho(S) = 2$. By Remark \ref{rho=2} we have 
$L_{|S} \eqv a C + b f$. Also $L_{|C}$ is very ample on $C \cong B$, therefore $d(B) \leq L.C = b - an$. Then 
$K_S + L_{|S} - f_1 - \ldots - f_{d(B)} \eqv aC + (b + n + 2g - 2 - d(B))f$. Hence $b + n + 2g - 2 - d(B) \geq n + 2g - 2 +
an \geq an + 2g - 1$, therefore $H^1(S, K_S + L_{|S} - f_1 - \ldots - f_{d(B)}) = 0$ by Lemma \ref{van1}(i). The vanishing
can be easily proved also in other cases, for example in the cases in \ref{uno}, \ref{due}, \ref{tre}, \ref{quattro} (either
with $B \cong \PP^1$ or with $\E$ sufficiently ample).}  
\end{rem}

To apply Corollary \ref{corestenell} to elliptic surfaces we need to verify the hypothesis on the
restriction map of Picard groups. It is here that the hypothesis on the rank becomes important. 

\begin{prop} 
\label{pic}
Let $X$ be a projective irreducible l.c.i.\ threefold, let $L$ be an ample line bundle on $X$ and let $S \in
|L|$ be a smooth surface with $\rho(S) = 2$ and $\kappa(S) = 1$. Then the restriction maps $\Pic X \to \Pic S$
and $N^1(X) \to N^1(S)$ are isomorphisms.
\end{prop}

\begin{proof}
Since $S$ is smooth we have that $\dim \Sing(X) \leq 0$ and $X$ is normal, whence $r_S : \Pic X \to \Pic S$ is
injective with torsion free cokernel by Lefschetz's theorem \cite[Cor.2.3.4]{bs}. Moreover the adjunction formula $K_S = (K_X +
L)_{|S}$ holds (see for example \cite[Prop.2.3 and 2.4]{ak}). Since $\kappa(S) = 1$ we have that
$L_{|S}$ and ${K_X}_{|S}$ are numerically independent and since $\rho(S) = 2$ we get that for any line bundle $A
\in \Pic S$ there are integers $a, u, v$ with $a \geq 1$ such that $aA \eqv u L_{|S} + v {K_X}_{|S}$. Therefore
we can write $aA \sim u L_{|S} + v {K_X}_{|S} + D$ with $D \eqv 0$. By \cite[Thm.4.6]{kl} there is an integer
$m \geq 1$ such that $mD \in \Pic^0 S$ and by \cite[Thm.2.3.1 and Thm.2.2.4]{bs} we have that $\Pic^0 X \to \Pic^0
S$ is an isomorphism, whence $mD \in \Im r_S$. Therefore also $maA \in \Im r_S$, whence $A \in \Im r_S$, since
$\Coker r_S$ is torsion free. 

The map $N^1(X) \to N^1(S)$ is now clearly surjective. To see its injectivity let $M \in \Pic X$ such that
$M_{|S} \eqv 0$. As above there is an integer $m \geq 1$ and a line bundle $N \in \Pic^0 X$ such that $mM_{|S}
\cong N_{|S}$, whence $mM \cong N$, therefore $M \eqv 0$. 
\end{proof}

The previous two results, together with some facts already present in the literature, allow us to give our best
version of an extension theorem for elliptic fibrations.

\begin{thm} 
\label{extell}
Let $X$ be a projective irreducible threefold, let $L$ be an ample line bundle on $X$ and let $S \in
|L|$ be a smooth surface having an elliptic fibration $\pi: S \to B$ with general fiber $f$. Then $\pi$
extends to a morphism ${\overline \pi} : X \to B$ if one of the following is satisfied:
\begin{itemize}
\item[(i)] $X$ is l.c.i.\ with rational singularities and $g(B) > 0$.
\item[(ii)] $X$ is l.c.i.\ with $\QQ$-factorial terminal singularities, not a cone over $S$, $L$ is very ample, $B \cong \PP^1$,
$\kappa(S) = 1$ and $N^1(S) \cong \ZZ [C] \oplus \ZZ [f]$ for some divisor $C$ such that $C.f \geq 1$. 
\item[(iii)] $X$ is l.c.i., $L$ is very ample, $\rho(S) = 2$, $\kappa(S) = 1$ and $H^1(S, K_S + L_{|S} - f_1 - \ldots - f_{d(B)}) =
0$ for every set of smooth distinct fibers $f_i$'s.
\item[(iv)] $X$ is smooth, $B \cong \PP^1$, $N^1(S) \cong \ZZ [C] \oplus \ZZ [f]$ for some divisor
$C$ and $L_{|S} \eqv aC + bf$ for any $a, b$ with $b > \frac{a C.f}{2} + 1 + \frac{1}{2a C.f} -
\frac{a C^2}{2 C.f}$. 
\end{itemize} 
\end{thm}

\begin{proof} Observe that (i) follows by \cite[Thm.5.2.3]{bs} since, as in the beginning of the proof of
Proposition \ref{pic}, $X$ is normal. Now (iii) is a consequence of Corollary \ref{corestenell} and Proposition
\ref{pic}. To see (iv) note that $a C.f = L_{|S}.f \geq 1$. Therefore the assumed inequality on $b$ is equivalent to 
$(f.L_{|S} + 1)^2 < L_{|S}^2$, whence $\pi$ extends by \cite[Thm.1.4]{pa}.

To prove (ii) we use first some adjunction theory to show that $K_X + L$ is nef.

By Proposition \ref{pic} we have that $\Pic X \to \Pic S$ and $N^1(X) \to N^1(S)$ are isomorphisms, so that $\rho(X) =
2$. Since $K_S \eqv ef$ for some $e \geq 1$, we have that $X$ does not admit a surjective morphism onto a variety $Y$ in
such a way that a general fiber $F$ intersects $S$ in a curve $\gamma$ with $p_a(\gamma) = 0$, for otherwise we would
have that $-2 = \gamma.(\gamma + K_S) = e \gamma.f \geq 0$. 

Now suppose that $K_X + L$ is not nef, so that $K_X$ is not nef and let $\tau$ be the nefvalue of $(X, L)$ (see Definition 
\ref{nefval}). Then $\tau > 1$ and by Proposition \ref{7.2.2} we deduce that $\tau \leq 3$ and that $K_X + 3L$ is nef and big. By
Theorem \ref{7.2.34} it follows that $K_X + 3L$ is ample and that $\tau \leq 2$. 

To go further note that we cannot have $K_X \sim - 2 L$, for otherwise $K_S \sim - L_{|S}$, giving the contradiction
$\kappa(S) = - \infty$. Furthermore let us show that $X$ cannot have a surjective morphism $\psi : X \to Y$ with
connected fibers onto a variety $Y$ of dimension $m = 1, 2$ in such a way that $K_X + 2 L \sim \psi^{\ast} \L$, for
some ample line bundle $\L \in \Pic Y$ such that $\psi$ is defined by the linear system $|q(K_X + 2 L)|$ for $q >> 0$.

In fact if such a $\psi$ exists then $\psi_{|S}$ must be surjective and if $m = 1$ we get the contradiction $0 < (K_S +
L_{|S})^2 = (K_X + 2 L)_{|S}^2 = (\psi^{\ast} \L)_{|S}^2 = (\psi_{|S}^{\ast} \L)^2 = 0$. On the other hand suppose that 
$m = 2$ and let $F_{\eta}$ be a general fiber of $\psi$. Then $q(K_X + 2 L)_{|F_{\eta}} \sim 0$ and by
\cite[Lemma3.3.2]{bs} we get $(K_X + 2 L)_{|F_{\eta}} \sim 0$, whence $K_{F_{\eta}} + 2 L_{|F_{\eta}} \sim 0$ and
therefore $(F_{\eta}, L_{|F_{\eta}}) \cong (\PP^1, \O_{\PP^1}(1))$. Now let $F$ be any fiber of $\psi$. Then $F$ is a line and if $F
\subset S$ we get the contradiction $0 \leq F.K_S = F.(K_X +  L)_{|S} = F.(K_X +  L) = - L.F < 0$. Now consider the Fano variety
$\F(X)$ of lines contained in $X$, which is at least $2$-dimensional as the fibers of $\psi$ are lines. Let $\Sigma_Y
\subseteq \F(X)$ be the irreducible subvariety of dimension two defined by the fibers of $\psi$. We have proved that the
hyperplane $H \subset \PP^{N+1} = \PP H^0(X, L)$ defining $S$ does not contain any line in $\Sigma_Y$. As is well known,
this implies that there is a point $P
\in \PP^{N+1}$ belonging to all lines in $\Sigma_Y$, therefore $X$ is a cone with vertex $P$. On the other hand, as $S$ 
is smooth, we get that $P \not\in S$, therefore $X$ is a cone over $S \subset \PP^N$, a contradiction.

Now by Theorem \ref{7.3.2} we have that either $K_X + 2L$ is ample or there exists a surjective birational morphism
$\phi : X \to X'$ onto a normal projective variety $X'$ with the following two properties: a) $\phi$ is the simultaneous
contraction to distinct smooth points $p_i' \in X'$ of some divisors $E_i \subset X - \Sing(X)$ such that $E_i \cong
\PP^2$, $\O_{E_i}(E_i) \cong \O_{\PP^2}(-1)$ and $L_{|E_i} \cong \O_{\PP^2}(1)$; b) If $L' := (\phi_{\ast} L)^{\ast
\ast}$ then $K_{X'} + 2 L'$ is ample. 

We now claim that if such a $\phi$ exists then it is an isomorphism. 

As a matter of fact suppose that $\phi$ is not an isomorphism, so that there is at least one contracted divisor $E_i$. 
Let $S' = \phi(S) \in |L'|$. As $S \cap \Sing(X) = \emptyset$ we have that $S'$ is certainly smooth outside all $p_i'$'s.
On the other hand ${E_i}_{|S}^2 = E_i.E_i.L = {E_i}_{|E_i}.L_{|E_i} = -1$, so that $p_i'$ is also a smooth point of $S'$.
Therefore $S'$ is birational, but not isomorphic to $S$, whence $\kappa(S') = 1$ and therefore $\rho(S') \geq 2$. But this
gives the contradiction $2 = \rho(S) \geq \rho(S') + 1 \geq 3$.

Hence, in any case, $K_X + 2L$ is ample, whence $1 < \tau < 2$ by \cite[Lemma1.5.5]{bs}. In this context the
first reduction $(\widetilde{X}, \widetilde{L})$ of $(X, L)$ is defined (see Section \ref{backgr}) and the above
discussion shows that in fact $(X, L) \cong (\widetilde{X}, \widetilde{L})$. Finally by Theorem \ref{7.3.4} we get that
$X$ has a surjective morphism onto a variety $Y$ in such a way that a general fiber $F$ intersects $S$ in a curve
$\gamma$ with $p_a(\gamma) = 0$, case already excluded.

This proves that $K_X + L$ is nef and the base-point free theorem (\cite[Thm.3.3]{km}) gives that 
for $q >> 0$ the line bundle $q(K_X + L)$ is base-point free. By \cite[Lemma1.1.3]{bs} there exists a surjective morphism
$p : X \to Y$ with connected fibers onto a normal projective variety $Y$ in such a way that $K_X + L \sim p^{\ast} \L$,
for some ample line bundle $\L \in \Pic Y$ and $p$ is defined by the linear system $|q(K_X + L)|$ for $q >> 0$. 

We now claim that $Y$ is a smooth curve. In fact let $F$ be a general fiber of $p$, so that $q(K_X + L)_{|F} \sim 0$,
whence, as above, $- K_F \sim L_{|F}$ by \cite[Lemma3.3.2]{bs}. If $Y$ is a point we find that $- K_X \sim L$, so that
$K_S \sim 0$, a contradiction. On the other hand we cannot have $\dim Y \geq 2$, for otherwise picking two general
divisors $D_1, D_2 \in |q(K_X + L)|$ we have $\dim D_1 \cap D_2 = 1$, giving the contradiction
\[ 0 = q^2 K_S^2 = q^2 (K_X + L)^2.L = D_1.D_2.L > 0. \]

Therefore $Y$ is a smooth curve and the general fiber $F$ of $p$ is a smooth connected surface. The 
restriction $p_{|S} : S \to Y$ is clearly surjective and a general fiber $f = F \cap S$ is therefore a smooth
connected curve with $f \in |-K_F|$, that is an elliptic curve. Therefore $p_{|S} : S \to Y$ is an elliptic
fibration and, as the latter is unique, we deduce that $Y = B, p_{|S} = \pi$ and that $\pi$ extends to
$X$, as required.
\end{proof}

\begin{rem}
{\rm The proof of the extension of $\pi$ under hypothesis (ii) is inspired by \cite{lm}. In case $X$ is smooth it has been
proved first by Ionescu \cite{io} (see also \cite{ds}, \cite{fds}).}
\end{rem}

The above theorem gives that most threefolds studied in this article, provided that they have the appropriate singularities,
are Mori fiber spaces.

\begin{prop} 
\label{mfs}
Let $X$ be a projective irreducible $\QQ$-factorial terminal l.c.i.\ threefold, let $L$ be a very ample line bundle on
$X$ and let $S \in |L|$ be a smooth surface having an elliptic fibration $\pi: S \to B$ with general fiber $f$. Suppose that
$\kappa(S) = 1$, $N^1(S) \cong \ZZ [C] \oplus \ZZ [f]$ for some divisor $C$ such that $C.f \geq 1$ and either $g(B) > 0$ or
$B \cong \PP^1$ and $X$ is not a cone over $S$. Then $\pi$ extends to a morphism ${\overline \pi} : X \to B$ and $X$ is a Mori fiber
space.
\end{prop}

\begin{proof}
Since terminal singularities are rational \cite[Thm.1]{elk}, by Theorem \ref{extell} $\pi$ extends to ${\overline \pi} : X
\to B$ with general fiber $F$ and by Proposition \ref{pic} we have that $N^1(X) \cong \ZZ [E] \oplus \ZZ [F]$ for a
divisor $E$ on $X$ such that $C = E_{|S}$. Hence we can write $K_X \eqv u E + v F$ and $L \eqv a E + b F$, for some
integers $u, v, a$ and $b$. Since $L$ is very ample we get $0 < L_{|S}.f = a C.f = a E_{|S}.F_{|S} = a E.F.(aE + bF)
= a^2 E^2.F$, so that $a \neq 0$ and $E^2.F \neq 0$. Now $0 = K_S.f = ((u + a)E + (v + b)F).F.L = a (u + a) E^2.F$
implies that $u = -a$ and therefore $- {K_X}_{|F} \eqv a E_{|F} \eqv L_{|F}$ is ample, whence $X$ is a Mori fiber
space. 
\end{proof}

We now study the fibers of the extended morphism.

\begin{thm}
\label{fibra}
Let $S \subset \PP^N$ be a smooth surface having an elliptic fibration $\pi: S \to B$ with general fiber $f$. Let
$X \subset \PP^{N+1}$ be a l.c.i.\ extension of $S$ given by a very ample divisor $L$ on $X$. Suppose
that $\pi : S \to B$ extends to a morphism $\overline{\pi} : X \to B$ and that $N^1(S) \cong \ZZ [C] \oplus \ZZ [f]$ 
for some divisor $C$ such that $C.f \geq 1$. Set $L_{|S} \eqv aC + bf$. Then
\begin{equation}
\label{possibili}
(a, C.f) \in \{ (1, 3), (1, 4), (1, 5), (1, 6), (1, 7), (1, 8), (1, 9), (2, 4), (3, 3) \}
\end{equation}
and for the general fiber $F$ of $\overline{\pi}$, we have $L_{|F} \cong -K_F$ and setting $d = K_F^2$, $F$
is one of the following:
\begin{itemize}
\item[(i)] $F \cong \PP^2$.
\item[(ii)] $F \cong \PP^1 \times \PP^1$.
\item[(iii)] $F$ is isomorphic to the blow-up of $9 - d$ distinct points in $\PP^2$ (no three collinear, no six on
a conic) for $3 \leq d \leq 8$.
\end{itemize}
Moreover if, in addition, $X$ is locally factorial, then $3 \leq d \leq 6$ in {\rm (iii)} and $(a, C.f) \not\in \{(1,
7), (1, 8), (1, 9) \}$ in \eqref{possibili}. 
\end{thm}

\begin{proof} We have $F_{|F} \eqv 0$, whence $\rho(X) \geq 2$. On the other hand, by \cite[Cor.2.3.4]{bs}, 
$\Pic X \to \Pic S$ is injective with torsion free cokernel. Therefore also $N^1(X) \to N^1(S)$ is injective with
torsion free cokernel, whence $N^1(X) \to N^1(S)$ is an isomorphism, since $\rho(S) = 2$. Let $[E]$ and $[F]$ be the
generators of $N^1(X)$, restricting respectively to $[C]$ and $[f]$ on $S$, so that $L \eqv aE + bF$ and $a C.f = L_{|S}.f
\geq 3$, giving $a \geq 1$.

We now claim that a general $F \subset \PP^{N+1}$ is a smooth Del Pezzo surface, that is $\O_F(1) \cong \O_F(-K_F)$. 

A general $F$ is certainly smooth by Bertini's theorem. Moreover $F \cap S = f$ is a smooth irreducible elliptic
curve, whence also $F$ is irreducible and it follows from \cite[Thm.8.9.3]{bs} that $F \subset \PP^{N+1}$ embedded by
$L_{|F}$ is either a Del Pezzo surface or $(F, L_{|F}) \cong (\PP \E, \xi)$, where $\E$ is a rank two vector bundle on
an elliptic curve and $\xi$ is the tautological line bundle.

Suppose we are in the latter case and let $\gamma$ be a fiber of $\PP \E$, so that $K_F \eqv -2 \xi + e \gamma$ for some
$e \in \ZZ$. We have $\xi = L_{|F} \eqv aE_{|F}$, whence $1 = \xi.\gamma = a E_{|F}.\gamma$ gives $a = 1$ and $\xi \eqv
E_{|F}$. Note that $X$ is Gorenstein, whence $K_X$ is a line bundle and we can write $K_X \eqv uE + vF$. Now 
\[ 0  = K_S.f = (K_X + L)_{|S}.f = (u + 1) C.f \]
therefore $u = -1$, giving the contradiction $-2 = K_F.\gamma = (K_X)_{|F}.\gamma = - E_{|F}.\gamma = - \xi.\gamma
= - 1$. 

This proves the claim and henceforth $F$ is a smooth Del Pezzo surface with $- K_F \sim L_{|F} \eqv aE_{|F}$. Now $d =
K_F^2 = L_{|F}^2 = L^2.F = L_{|S}.f = a C.f$ and by \cite[Thm.8]{nag} we have that $3 \leq d \leq 9$. 

By \cite[Thm.8]{nag} we deduce that either $a = 3$ and $F \cong \PP^2$, so that $C.f = 3$ or $a = 2$ and $F \cong \PP^1
\times \PP^1$, so that $C.f = 4$ or $a = 1$ and $3 \leq C.f \leq 9$. Thus \eqref{possibili} is proved.

Moreover, again by \cite[Thm.8]{nag}, when $a = 1$, we have that either $L_{|F}$ is divisible, leading to $C.f = 9$, $F
\cong\PP^2$ and $C.f = 8$, $F \cong \PP^1 \times \PP^1$ or $L_{|F}$ is not divisible. 

Suppose we are in the latter case. By \cite[Thm.8]{nag} we get $3 \leq d \leq 8$ and
$F$ is the anticanonical embedding of the blow-up of $\PP^2$ in $9 - d$ points. Note that, as $F$ is smooth, the
blown-up points cannot be infinitely near, there is no line containing three of them and there is no conic containing
six of them.

This proves that $F$ is as in (i), (ii) or (iii).

For the rest of the proof suppose that $X$ is also locally factorial. 

If $F$ is as in (iii), there is a nonempty open subset $U \subset B$ such that each fiber over $U$ is as in (iii). This gives
rise to an irreducible curve $B' \subset \F(X)$, the Fano variety of lines contained in $X$. Let $T$ be the union of all the
lines in $B'$. By \cite[Exa.6.19 and Prop.6.13]{hr} we have that $T$ is a Weil divisor on $X$ and by \cite[Prop.II.6.11]{ht}
it follows that $T$ is also a Cartier divisor on $X$. Hence $T \eqv rE + sF$ for some integers $r$ and $s$. Therefore, using
the convention that $\binom{m}{n} = 0$ if $m < n$, we have
\[ 9 - d + \binom{9 - d}{2} + \binom{9 - d}{5} = \deg T_{|F} = L_{|F}.T_{|F} =  L.T.F =  T_{|S}.F_{|S} = 
(rC + sf).f = rd, \] 
giving $3 \leq d \leq 6$. 

Finally we exclude the cases $(a, C.f) \in \{(1, 7), (1, 8), (1, 9) \}$ in \eqref{possibili}. If $(a, C.f) = (1, 7)$ we
have $d = a C.f = 7$, contradicting what we have just proved. Now assume that either $(a, C.f) = (1, 9)$, so that 
$F \cong \PP^2$ or $(a, C.f) = (1, 8)$ and $F \cong \PP^1 \times \PP^1$. Let $k(B)$ be the quotient field of $B$ and
$\overline{k(B)}$ its algebraic closure, so that the base changes of $X$, $X_{k(B)}$ to $k(B)$ and $X_{\overline{k(B)}}$
to $\overline{k(B)}$ are defined. Now either $X_{\overline{k(B)}} \cong \PP_{\overline{k(B)}}^2$ or $X_{\overline{k(B)}}
\cong \PP_{\overline{k(B)}}^1 \times \PP_{\overline{k(B)}}^1$, and, as in \cite[Proof of 3.5.2, page 162]{mor}, we have
$\Pic X_{k(B)} \cong (\Pic X_{\overline{k(B)}})^G$, where $G = \Gal(\overline{k(B)}/k(B))$ is the Galois group. This
implies that the canonical bundle $K_{X_{k(B)}}$ is $r$ divisible for $r = 3$ in the case $F \cong \PP^2$ and $r = 2$ in
the case $F \cong \PP^1 \times \PP^1$. Therefore we can find a nonempty open subset $U \subset B$ and a line bundle $\L$
on $\overline{\pi}^{-1}(U)$ such that $\L_{|F_u} \cong \O_{\PP^2}(1)$ if $F_u \cong \PP^2$ and $\L_{|F_u} \cong \O_{\PP^1
\times \PP^1}(1, 1)$ if $F_u \cong \PP^1 \times \PP^1$ on every fiber $F_u$ over $U$. Since $X$ is locally factorial, $\L$
extends to a line bundle $\overline{\L}$ having the same restriction property as $\L$ on a general fiber $F$. But now
$\overline{\L} \eqv \alpha E + \beta F$, for some $\alpha, \beta \in \ZZ$ and therefore 
\[\L_{|F} \eqv \overline{\L}_{|F} \eqv \alpha E_{|F} \eqv \alpha L_{|F} \eqv - \alpha K_F \]
which easily gives a contradiction.
\end{proof}

With this baggage of results we are now ready to prove our main theorems.

\renewcommand{\proofname}{Proof of Theorem {\rm \ref{nonextell}}}  
\begin{proof}
Immediate consequence of Theorems \ref{extell} and \ref{fibra}.
\end{proof}
\renewcommand{\proofname}{Proof}

\renewcommand{\proofname}{Proof of Corollary {\rm \ref{nonextweier}}}  
\begin{proof}
If $(g, n) \neq (0, 2)$ the Corollary follows from Proposition \ref{invariants}(iv), Remark \ref{rho=2}, Theorem
\ref{nonextell}(iii) and Remark \ref{van}. Now suppose that $(g, n) = (0, 2)$, so that $S$ is a $K3$ surface
by Proposition \ref{invariants}(iii). In particular any line bundle numerically equivalent to $0$ on $S$ is
$\O_S$. Since $S$ is smooth we have that $X$ is normal and $\dim \Sing(X) \leq 0$, whence $r_S : \Pic X \to \Pic S
\cong N^1(S)$ is injective with torsion free cokernel by \cite[Cor.2.3.4]{bs}. Let $A \in \Pic X$ be such that $A
\eqv 0$. Then $A_{|S} \eqv 0$, whence $A_{|S} \sim 0$ and therefore $A \sim 0$ by the injectivity of $r_S$. Hence
also $\Pic X \cong N^1(X)$. Therefore $\rho(X) \leq 2$. Now if $\rho(X) = 2$ then necessarily $\Pic X = N^1(X)
\cong N^1(S) = \Pic S$, therefore the morphism $\pi : S \to B$ extends to a morphism $\overline{\pi} : X \to B$
by Corollary \ref{corestenell} and Remark \ref{van}. Now we get a contradiction by Theorem \ref{fibra}.
 
Hence $\rho(X) = 1$. By \cite[Thm.5.3.1]{bs} we also get that $-K_X \in |\O_X(1) |$ and that $h^1(\O_X) =
h^2(\O_X) = 0$.
\end{proof}
\renewcommand{\proofname}{Proof}

\begin{rem}
{\rm In the case (ii) of Corollary \ref{nonextweier}, if $X$ is smooth, the surface $S$ must belong to some proper
closed subset of the linear system $|L|$: In fact, for a general $S' \in |L|$ we have $1 = h^2(\O_{S'}) >
h^2(\O_X) = 0$, whence $\Pic X \cong \Pic S'$ by a theorem of Moishezon \cite[Thm.7.5]{moi}.}
\end{rem}

We now consider the extendability of $K3$ Weierstrass fibrations.

\begin{lemma}
\label{>12}
Let $S \subset \PP^N$ be a smooth Weierstrass fibration with $\kappa(S) = 0$ and $\rho(S) = 2$.
Then the sectional genus $g(S)$ of $S$ satisfies $g(S) \geq 13$.
\end{lemma}

\begin{proof} As usual let $C$ be the section and $f$ be the general fiber. Recall that, by
Proposition \ref{invariants}(iii) we have $C^2 = -2$. If $H_S$ is the hyperplane bundle of
$S$ we have by Remark \ref{rho=2} that $H_S \sim aC + bf$ with $a \geq 3$ and $b \geq 2a + 1$ 
by Lemma \ref{van1}(vi), whence
\[ g(S) = \frac{1}{2} H_S^2 + 1 = - a^2 +  ab + 1 \geq a(a + 1) + 1 \geq 13. \]
\end{proof}

\renewcommand{\proofname}{Proof of Corollary {\rm \ref{nonextweierk3}}}  
\begin{proof}
Let $X \subset \PP^{N+1}$ be an extension of $S \subset \PP^N$. If $X$ is normal we have by Lemma \ref{>12} and
\cite[Cor.1.6]{pr}, that $13 \leq g(S) = a(b - a) + 1 \leq 37$. Using Lemma \ref{van1}(vi) we get (iii). 

To see (i) and (ii) suppose that $X$ is l.c.i.. By Corollary \ref{nonextweier} we have that $X$ is an anticanonically
embedded Fano threefold with Picard number one. Set $H = - K_X$. Moreover, as in the proof of Corollary
\ref{nonextweier}, we know that $\Pic X \to \Pic S$ is injective with torsion free cokernel.

To show (ii) we will exclude from the list in (iii) the five cases 
\[(a, b, g(S)) \in \{(3, 9, 19), (3, 12, 28), (3, 15, 37), (4, 10, 25), (4, 12, 33) \}.\] 
As a matter of fact in the above cases we have that $H_S$ is $r$-divisible in $\Pic S$ with $r = 3$ in the first three
cases, $r = 2$ in the fourth case and $r = 4$ in the fifth case. Hence $H \sim r \Delta$ for some ample $\Delta \in \Pic X$.
By the generalized Kobayashi-Ochiai theorem \cite[Thm.3.1.6]{bs} we deduce that $(X, \Delta) \cong (\PP^3, \O_{\PP^3}(1))$
when $r = 4$, while $(X, \Delta) \cong (Q, \O_Q(1))$, where $Q \subset \PP^4$ is a quadric when $r = 3$. If
$r = 2$ we know that $\Delta_{|S} \sim 2C + 5 f$ and $(X, \Delta)$ is a Del Pezzo variety. In the above case we claim that
$H^i(t \Delta) = 0$ for $0 < i < 3$ and for every $t \in \ZZ$. To see this consider the exact sequence
\begin{equation}
\label{seqk3}
0 \hpil (t - 2) \Delta \hpil t \Delta \hpil t \Delta_{|S} \hpil 0.
\end{equation}
By Lemma \ref{van1}(i), Serre duality and the fact that $S$ is a K3 surface we know that $H^1(t \Delta_{|S}) = H^1( 2tC +
5t f) = 0$ for every $t \in \ZZ$ and $H^2(t \Delta_{|S}) = H^2( 2tC + 5t f) = 0$ for every $t \geq 1$. Therefore we
deduce from \eqref{seqk3} that (a) $h^1(t \Delta) \leq h^1((t - 2) \Delta)$ for every $t \in \ZZ$ and (b) $h^2((t - 2)
\Delta) = h^2(t \Delta)$ for every $t \geq 1$. Since $h^2(t \Delta) = 0$ for $t >> 0$, using (b), we find that $h^2(t
\Delta) = 0$ for every $t \geq -1$ and therefore $h^1(t \Delta) = 0$ for every $t \leq -1$. This, together with (a), gives
that
$h^1(t \Delta) = 0$ for every $t \in \ZZ$ and also $h^2(t \Delta) = 0$ for every $t \in \ZZ$ by Serre duality.
This proves the claim and it follows from \cite[Cor.1.5]{f} that, also in this case, $\Delta$ is very ample.

Hence, in all cases, $\Delta_{|S}$ is very ample and therefore, by Lemma \ref{van1}(vi), we get that $a \geq 3r \geq
6$, a contradiction. This proves (ii).

Finally, to prove (i), suppose that $X$ is also with terminal singularities. By \cite[Thm.11]{nam} we have that $X$ is
smoothable, whence a general deformation $X_{\eta}$ of $X$ is a smooth Fano threefold with Picard number one and with $-
K_{X_{\eta}}$ very ample. Moreover $X_{\eta}$ has genus $g = g(S) \geq 13$ by Lemma \ref{>12}. By \cite[Thm.6]{clm1} and
\cite[Thm.3.2]{clm2} or by \cite[Thm.4.2]{isk1}, \cite[Thm.6.1]{isk2} (together with \cite{s1, s2}), we find that
\[ g(S) = a(b - a) + 1 = 13, 17, 21, 28, 33 \]
and also that $- K_{X_{\eta}}$ is $2$-divisible in the first three cases. Therefore also $H = - K_X$ has the same
divisibility properties, whence so  does $aC + b f \eqv H_S = H_{|S}$. By Lemma \ref{van1}(vi) we know that $a \geq 3$
and $b \geq 2a + 1$, whence, using (ii), we get the two possibilities $(a, b, g(S)) \in \{(3, 7, 13), (4, 9, 21) \}$. As
we said above, $H_S$ is $2$-divisible, contradicting $H_S.C = 1$ in the case $(4, 9, 21)$ and $H_S.f = 3$ in the
case $(3, 7, 13)$. This proves (i). 
\end{proof}
\renewcommand{\proofname}{Proof} 

\section{Nonextendability of many embeddings of Weierstrass fibrations}

In \cite{glm} and \cite{klm} a new technique to deal with the extendability of a surface was introduced.
It is the purpose of this section to recall it and then use it to prove a nonextendability result (regardless of
the singularities) for many very ample line bundles on a Weierstrass fibration.

We first recall the definition and notation for multiplication maps and Gaussian maps.

\begin{notn}
Let $L, M$ be two line bundles on a smooth projective variety $X$. Given $V \subseteq H^0(L)$ we will denote by
$\mu_{V, M} : V \otimes H^0(M) \hpil H^0(L \otimes M)$ the multiplication map of sections, $\mu_{L, M}$ when $V
= H^0(L)$, by $R(L, M)$ the kernel of $\mu_{L, M}$ and by $\Phi_{L,M} : R(L, M) \hpil  H^0(\Omega^1_X \otimes L
\otimes M)$ the Gaussian map (that can be defined locally by $\Phi_{L,M}(s \otimes t) = sdt - tds$, see
\cite[1.1]{w}).
\end{notn}

Let us also recall, for the reader's convenience, a couple of results about Gaussian maps that will be used in the sequel.

\begin{prop} \cite[Prop.1.10]{w}
\label{1.10}
Let $L$ be a very ample line bundle on a smooth irreducible variety $X$ giving an embedding $X \subset \PP^r$ and let $M$ be another
line bundle. Then
\begin{itemize}
\item[(i)] the Gaussian map $\Phi_{L, M}$ is the restriction map $H^0(X, \Omega^1_{\PP^r} \otimes \O_X \otimes L \otimes M) \to H^0(X,
\Omega^1_X \otimes L \otimes M)$; 
\item[(ii)] $\Coker \Phi_{L, M} = \Ker \{ H^1(X, N_{X/\PP^r}^{\ast} \otimes L \otimes M) \to H^1(X, \Omega^1_{\PP^r} \otimes \O_X
\otimes L \otimes M) \}$; 
\item[(iii)] if $\mu_{L, M}$ is surjective and $H^1(M) = 0$ then $\Coker \Phi_{L, M} \cong H^1(X, N_{X/\PP^r}^{\ast} \otimes L
\otimes M)$. 
\end{itemize}
\end{prop}

\begin{thm} \cite[Thm.2]{bel}
\label{Thm.2}
Let $L$ be a line bundle on a smooth irreducible curve $C$ of genus $g$. If $\Cliff(C) \geq 2$ and $\deg L \geq 4g + 1 - 2
\Cliff(C)$ or if $\Cliff(C) \geq 3$ and $\deg L \geq 4g + 1 - 3 \Cliff(C)$ then $\Phi_{\omega_C, L}$ is surjective.
\end{thm}

The way these maps are used to prove nonextendability is explained in the following

\begin{prop} \cite[Cor.2.4]{klm}
\label{cor:anysurface}
Let $Y \subset \PP^r$ be a smooth irreducible surface which is either linearly normal or regular
(that is, $h^1(\O_Y) = 0$) and let $H$ be its hyperplane bundle. Assume there is a base-point
free and big line bundle $D_0$ on $Y$ with $H^1(H-D_0) = 0$ and such that the general element $D
\in |D_0|$ is not rational and satisfies
\begin{itemize}
\item[(i)] the Gaussian map $\Phi_{H_D, \omega_D}$ is surjective; 
\item[(ii)] the multiplication maps $\mu_{V_D, \omega_D}$ and $\mu_{V_D, \omega_{D}(D_0)}$ are
surjective, where 

$V_D := \Im \{H^0(Y, H-D_0) \khpil H^0(D, (H-D_0)_{|D})\}$.
\end{itemize}
Then $Y$ is nonextendable. 
\end{prop}

Now, on a Weierstrass fibration, we can translate the hypotheses of Proposition \ref{cor:anysurface} into
(essentially) purely numerical conditions.

\begin{prop}
\label{numerical}
Let $S \subset \PP^N$ be a smooth surface having a Weierstrass fibration $\pi: S \to B$ with general fiber $f$
and section $C$ and with $n = - C^2 \geq 1$ and $g = g(B)$. Suppose that either $g = 0$ or $S$ is linearly normal
and that the hyperplane bundle of $S$ is of type $H \eqv aC + bf$. 

Given integers $\alpha$ and $\beta$ let $D_0 = \alpha C + \beta f$. Set $\W = 0$ when $g = 0$ and $\W \in \Pic^0 B$
general (with respect to $D_0$) when $g \geq 1$. Let $D_{0, \W} = D_0 \otimes \pi^{\ast} \W$ and suppose that the
general $D \in |D_{0, \W}|$ is not hyperelliptic and satisfies
\begin{itemize}
\item[(a)] $10 \neq D_0.(D_0 + K_S) \geq 6$;
\item[(b)] $\alpha \geq 2$ and $\beta \geq \alpha n + 2g$;
\item[(c)] either $a - 2\alpha \geq 2$ and $b - 2\beta \geq (a - 2 \alpha) n + g - 1$ or $a - 2\alpha = 1$
and $b - 2\beta \geq g - 1$;
\item[(d)] $(H-D_0).D_0 \geq D_0.(D_0 + K_S) + 3$;
\item[(e)] If $D$ is not trigonal then $H.D_0 \geq 2D_0.(D_0 + K_S) + 1$;
\item[(f)] If $D_0.(D_0 + K_S) \geq 8$ and $D$ is trigonal then $H.D_0 \geq \frac{3}{2} D_0.(D_0 + K_S) + 10$;
\item[(g)] If $D_0.(D_0 + K_S) = 6$ then $H.D_0 \geq 17$.
\end{itemize}
Then $S$ is not extendable.
\end{prop}

\begin{proof} Since $D_{0, \W}^2 = D_0^2 = \alpha (2\beta - \alpha n)$, by (b) it follows that $D_{0, \W}^2 \geq
\alpha (\alpha n + 4g) > 0$. By Lemma \ref{van1}(iv) and again (b) we know that $|D_{0, \W}|$ is base-point free,
whence $D$ is smooth and irreducible by Bertini's theorems. By (a) we have that $6 \neq g(D) \geq 4$, in particular
$D$ is not isomorphic to a plane quintic.

By (c) and Lemma \ref{van1}(iii) we deduce that $H^1(H - 2D_{0, \W}) = 0$ and therefore that $V_D = H^0((H -
D_{0, \W})_{|D})$. Now (d) is just $(H - D_{0, \W}).D \geq 2g(D) + 1$, whence $H^1((H - D_{0, \W})_{|D}) = 0$ and
$(H - D_{0, \W})_{|D}$ is very ample. Therefore $\mu_{V_D, \omega_D} = \mu_{(H - D_{0, \W})_{|D}, \omega_D}$ is
surjective by \cite[Thm.1.6]{as} and the exact sequence
\[0 \hpil H - 2D_{0, \W} \hpil H - D_{0, \W} \hpil (H - D_{0, \W})_{|D} \hpil 0 \]
shows that also $H^1(H - D_{0, \W}) = 0$. 

The multiplication map $\mu_{V_D, \omega_{D}(D_{0, \W})} = \mu_{(H - D_{0, \W})_{|D}, \omega_{D}(D_{0, \W})}$ is
surjective by Green's $H^0$-lemma \cite[Thm.4.e.1]{gr}: In fact we just need to verify that  
\[ h^0((H - D_{0, \W})_{|D} - {D_{0, \W}}_{|D}) = h^1(\omega_D(D_{0, \W}) - (H - D_{0, \W})_{|D}) \leq h^0((H -
D_{0, \W})_{|D}) - 2, \]
which holds since $(H - D_{0, \W})_{|D}$ is very ample and ${D_{0, \W}}_{|D}$ is effective of degree at least $2$.

To apply Proposition \ref{cor:anysurface} it remains to check that the Gaussian map $\Phi_{H_D, \omega_D}$ is
surjective. Now if $g(D) = 4$ this follows by (g) and Proposition \ref{1.10} (or by \cite[Prop.2.9 ]{klgm}). If $g(D)
\geq 5$ and $D$ is trigonal, this follows by (f) and \cite[Cor.2.10]{klgm}. Finally if $D$ is not trigonal
this follows by (e) and Theorem \ref{Thm.2}.
\end{proof} 

We can also use the standard scroll containing any Weierstrass fibration to compute the cohomology of the normal bundle.

\begin{lemma} 
\label{nonextweieronscroll}
Let $S \subset \PP^N$ be a smooth surface having a Weierstrass fibration $\pi: S \to B$ with general fiber 
$f$ and section $C$. Set $n = - C^2$ and $g = g(B)$. Suppose that the hyperplane bundle of $S$ is of type $H_S
\eqv aC + bf$ and that $n \geq 1$. 

If $a = 3u$ for some $u \geq 2$ and $b \neq a + 1$ if $(n, g ) = (1, 0)$, then $H^1(T_S(-H_S)) =
0$, where $T_S$ is the tangent bundle of $S$.
\end{lemma}

\begin{proof} 
Let $\L$ be the fundamental line bundle of the fibration and let $\E = \pi_{\ast} \O_S(3C) \cong \O_B \oplus
\L^{-2} \oplus \L^{-3}$ and $Y = \PP \E$ be the threefold scroll with projection morphism $p: Y \to B$. By
\cite[III.1]{mi2} $S$ can be embedded as a divisor linearly equivalent to $3 \xi + 6 p^{\ast} \L$ in $Y$. As
$\xi_{|S} \eqv 3C$ we have that $H_S \sim A_{|S}$ for some line bundle $A \eqv u \xi + b F$ on $Y$.
Therefore there exists a line bundle $M \in \Pic^b B$ such that $A \sim u \xi + p^{\ast} M$. We will often
use the fact that, as $H_S$ is very ample, we have $H_S.C \geq 1$, whence
\begin{equation}
\label{limb}
b \geq an + 1.
\end{equation}
The goal will be to use the scroll $Y$ to compute cohomology on $S$.

\begin{claim} 
\label{claim1} 
With notation as above we have $H^0(N_{S / Y}(-H_S)) = 0$.
\end{claim}

\begin{proof} 
We have $N_{S / Y}(-H_S) \eqv (9 - 3u)C + (6n - b)f$, so that the required vanishing is obvious if $u \geq
4$, while it follows by pushing down to $B$ and using Lemma \ref{split}(ii) when $u = 2, 3$ since, in this case,
$b \geq 6n + 1$ by \eqref{limb}.
\end{proof}

\begin{claim} 
\label{claim2}  
With notation as above we have $H^1(p^{\ast}(- K_B)(- A)) = 0$ and $H^2(p^{\ast}(- K_B)(- A - S)) = 0$.
\end{claim}

\begin{proof} 
We have $p^{\ast}(-K_B)(-A) \sim - u \xi + p^{\ast}(- K_B - M)$. Since $R^i p_{\ast} (- u \xi) = 0$ 
for $i = 0, 1$ we get the first vanishing by the Leray spectral sequence. As for the second, by Serre duality,
we need to show that $H^1(p^{\ast}(K_B)(K_Y + A + S)) = 0$. Since $p^{\ast}(K_B)(K_Y + A + S) \sim u \xi +
p^{\ast}(2 K_B + \L + M)$ and $R^1 p_{\ast} (u \xi) = 0$, we deduce, again by the Leray spectral
sequence, that $H^1(Y, p^{\ast}(K_B)(K_Y + A + S)) \cong H^1(B, \Sym^u \E (2 K_B + \L + M)) = 0$ for degree
reasons (here we use \eqref{limb} and the hypothesis $b \neq a + 1$ if $(n, g ) = (1, 0)$).
\end{proof}

\begin{claim} 
\label{claim3} 
With notation as above we have $H^1(p^{\ast} \E^{\ast}(\xi - A)) = 0$ and $H^2(p^{\ast} \E^{\ast}(\xi - A -
S)) = 0$.
\end{claim}

\begin{proof} 
By Serre duality we have $H^1(p^{\ast} \E^{\ast}(\xi - A)) \cong H^2(p^{\ast} \E (K_Y - \xi + A))$
and $p^{\ast} \E (K_Y - \xi + A) \cong p^{\ast} (\E (K_B + M - 5 \L))((u - 4) \xi)$. Since $R^i p_{\ast} ((u
- 4) \xi) = 0$  for $i = 1, 2$ we get the first vanishing by the Leray spectral sequence. As for the second, 
by Serre duality, we need to show that $H^1(p^{\ast} \E (K_Y - \xi + A + S)) = 0$. Since $p^{\ast} \E (K_Y -
\xi + A + S) \cong p^{\ast} (\E (K_B + \L + M))((u - 1) \xi)$ and $R^1 p_{\ast} ((u - 1) \xi) = 0$, we deduce,
again by the Leray spectral sequence, that $H^1(Y, p^{\ast} \E (K_Y - \xi + A + S)) \cong H^1(B, \Sym^{u - 1}
\E \otimes \E (K_B + \L + M)) = 0$ for degree reasons (using \eqref{limb}).
\end{proof}

\begin{claim} 
\label{claim4} 
With notation as above we have $H^2(\O_Y(- A)) = 0$.
\end{claim}

\begin{proof} 
By Serre duality we have $H^2(\O_Y(- A)) \cong H^1(K_Y + A)$ and $K_Y + A \cong p^{\ast}
(K_B + M - 5 \L))((u - 3) \xi)$. Since $R^1 p_{\ast} ((u - 3) \xi) = 0$ and $p_{\ast} (- \xi) = 0$
we get, by the Leray spectral sequence, that $H^1(K_Y + A) = 0$ for $u = 2$ and $H^1(K_Y + A) \cong
H^1(B, \Sym^{u - 3} \E (K_B + M - 5 \L)) = 0$ for degree reasons (using \eqref{limb}) for $u \geq 3$.
\end{proof}

For the sequel we will let $T_{Y / B}$ be the relative tangent bundle and use the two standard exact sequences
\begin{equation}
\label{rel1}
0 \hpil T_{Y / B} \hpil T_Y \hpil p^{\ast}(- K_B) \hpil 0
\end{equation}
and
\begin{equation}
\label{rel2}
0 \hpil \O_Y \hpil p^{\ast} \E^{\ast}(\xi)  \hpil T_{Y / B} \hpil 0.
\end{equation}

\begin{claim} 
\label{claim5} 
With notation as above we have $H^1(T_{Y / B}(- A)) = 0$ and $H^2(T_{Y / B}(- A - S)) = 0$.
\end{claim}

\begin{proof} 
Tensoring \eqref{rel2} with $\O_Y(- A)$ we get the first vanishing by Claims \ref{claim3} and \ref{claim4}.
Tensoring \eqref{rel2} with $\O_Y(- A - S)$ we get a map 
\[ \varphi : H^3(\O_Y(- A - S)) \to H^3(p^{\ast} \E^{\ast}(\xi - A - S)). \]
By Claim \ref{claim3} and Serre duality we find that
\[ H^2(T_{Y / B}(- A - S)) \cong \Ker \varphi \cong \Coker \varphi^{\ast} \]
and we need to prove that $\varphi^{\ast} : H^0(p^{\ast} \E(K_Y - \xi + A + S)) \to H^0(K_Y + A + S)$ is
surjective. Pushing down to $B$ we see that this is equivalent to the surjectivity of
the natural multiplication map $H^0(B, \Sym^{u - 1} \E \otimes \E (K_B + \L + M)) \to H^0(B, \Sym^u \E (K_B +
\L + M))$. Now the latter is surjective since, as $\E$ is split, $\Sym^u \E$ is a direct summand of $\Sym^{u -
1} \E \otimes \E$ and the map is given by projection.
\end{proof}

\noindent {\it Conclusion of the proof of Lemma {\rm \ref{nonextweieronscroll}}}.
Tensoring \eqref{rel1} with $\O_Y(- A - S)$ we see by Claims \ref{claim5} and \ref{claim2} that 
\[ H^2(T_Y (- A - S)) = 0.  \]
Tensoring \eqref{rel1} with $\O_Y(- A)$ we see by Claims \ref{claim5} and \ref{claim2} that 
\[ H^1(T_Y (- A)) = 0.  \]
Now from the exact sequence
\[ 0 \hpil T_Y (- A - S) \hpil T_Y (- A) \hpil {T_Y}_{|S}(- H_S) \hpil 0 \]
we deduce that $H^1({T_Y}_{|S}(- H_S)) = 0$. Finally Claim \ref{claim1} and the exact sequence
\[ 0 \hpil T_S(- H_S) \hpil {T_Y}_{|S}(- H_S) \hpil N_{S / Y}(-H_S) \hpil 0 \]
prove the Lemma.
\end{proof} 

We are now ready to prove nonextendability.

\renewcommand{\proofname}{Proof of Theorem {\rm \ref{nonextweier2}}}  
\begin{proof}
We first apply Proposition \ref{numerical} with $D_0 = 3C + 3nf$ if $g = 0$ and $D_0 = 2C + (2n + 2g)f$ if $g \geq
1$. The fact that $D$ is not hyperelliptic is a consequence of Lemma \ref{nonhyper}. When $g = 0$ we know that $D$
is trigonal, while if $g \geq 1$, we have, again by Lemma \ref{nonhyper}, that $D$ is not trigonal. A
straightforward calculation now proves that $S$ is not extendable if any of the following conditions is satisfied:
\begin{eqnarray}
\label{prima}
\bullet & \ \ & g = 0 \mbox{ and } a \geq 7, \\ \nonumber \bullet & & g \geq 1, S \mbox{ is linearly normal and either } a
\geq 7, b \geq an + 5g - 1 \mbox{ or, } \\ \nonumber & & a = 6, b \geq \max \{ 6n + 5g - 1, 6n + 6g - 3 \} \mbox{ or
} a = 5, b \geq 6n + 7g - 3.
\end{eqnarray}
On the other hand we know by Lemma \ref{nonextweieronscroll} that if 
\begin{eqnarray}
\label{seconda}
\bullet \ \ a = 3u \mbox{ for some } u \geq 2 \mbox{ and } b \neq a + 1 \mbox{ if } (n, g ) = (1, 0),
\end{eqnarray}
then $H^1(T_S(- 1)) = 0$. Now the Euler sequence of $S \subset \PP^N$ implies that $h^0({T_{\PP^N}}_{|S}(- 1)) = N + 1$ and
therefore also $h^0({N_{S / \PP^N}}(- 1)) = N + 1$. We deduce by Zak's theorem \cite[page 277]{z} (see also \cite[Thm.0.1]{lv})
that $S$ is not extendable under condition \eqref{seconda}. Putting this together with \eqref{prima}, the theorem is proved.
\end{proof}
\renewcommand{\proofname}{Proof}

\section{Examples of extendable elliptic surfaces}
\label{eoes}

In this section we will exhibit some simple examples of smoothly extendable elliptic surfaces $S$ with $\rho(S) = 2$ and some 
examples of smoothly extendable Weierstrass fibrations $S$ with $\rho(S) = 3, 4$.

Let $B$ be a smooth curve and let $\E$ be a very ample vector bundle on $B$. We denote by $\xi$ the 
tautological line bundle on $\PP \E$ and by $F$ a fiber of $\pi : \PP \E \to B$. We recall that $N^1(\PP \E) \cong
\ZZ[\xi] \oplus \ZZ[F]$. Given a surface $S \subset \PP \E$ we denote by $C = \xi_{|S}$ and by $f = F_{|S}$.

\begin{exa} 
\label{uno}
Suppose $\E$ has rank $3$, let $X = \PP \E$, let $L \cong \O_X(3 \xi)$ and consider the embedding $X
\subset \PP H^0(L)  = \PP^{N+1}$. Let $S \in |L|$ be a general hyperplane section (in the countable Zariski
topology). Then $\pi_{|S} : S \to B$ is an elliptic fibration and $N^1(S) \cong \ZZ[C] \oplus
\ZZ[f]$ by \cite[Thm.7.5]{moi} since $p_g(S) = h^2(\O_S) > h^2(\O_X) = 0$. Here $L_{|S} \eqv 3C$, $C.f = 3$ and
$L_{|S}.f = 9$. With the notation of Theorem \ref{nonextell}, this is the case $(a, C.f) = (3, 3)$. 
\end{exa}

\begin{exa} 
\label{due}
Suppose $\E$ has rank $4$, let $Y = \PP \E$ and consider the embedding $Y \subset \PP H^0(\xi) = \PP^{N+1}$. Take
a general divisor $X \in |3 \xi|$. Then, by Gherardelli-Lefschetz's theorem \cite{ghe}, \cite[Cor.2.3.4]{bs}, we
have $N^1(X) \cong \ZZ[\xi_{|X}] \oplus \ZZ[F_{|X}]$. For the embedding $X \subset Y \subset \PP^{N+1}$,
the general hyperplane section (in the countable Zariski topology) $S \subset \PP^N$ of $X$ is an extendable
elliptic surface over $B$ and, as in Example {\rm \ref{uno}}, we have $N^1(S) \cong \ZZ[C] \oplus
\ZZ[f]$ by \cite[Thm.7.5]{moi}. If $L = \O_X(1)$ we get $L_{|S} \eqv C$, $C.f = 3$ and $L_{|S}.f = 3$. With the
notation of Theorem \ref{nonextell}, this is the case $(a, C.f) = (1, 3)$.
\end{exa}

\begin{exa} 
\label{tre}
Suppose $\E$ has rank $4$, let $Y = \PP \E$ and consider the embedding $Y \subset \PP H^0(2 \xi) = \PP^{N + 2}$.
Take a general hyperplane section $X = Y \cap H \subset \PP^{N+1}$. As in Example {\rm \ref{due}} we have
$N^1(X) \cong \ZZ[\xi_{|X}] \oplus \ZZ[F_{|X}]$. For the embedding $X \subset Y \subset \PP^{N+1}$,
the general hyperplane section (in the countable Zariski topology) $S \subset \PP^N$ of $X$ is an extendable
elliptic surface over $B$ and, as in Example {\rm \ref{uno}}, we have $N^1(S) \cong \ZZ[C] \oplus
\ZZ[f]$ by \cite[Thm.7.5]{moi}. If $L = \O_X(1)$ we get $L_{|S} \eqv 2C$, $C.f = 4$ and $L_{|S}.f = 8$. With the
notation of Theorem \ref{nonextell}, this is the case $(a, C.f) = (2, 4)$.
\end{exa}

\begin{exa} 
\label{quattro}
Suppose $\E$ has rank $5$, let $Y = \PP \E$ and consider the embedding $Y \subset \PP H^0(\xi) = \PP^{N + 1}$.
Take two general quadrics $Q_1, Q_2$ and let $X = Y \cap Q_1 \cap Q_2$. As in Example {\rm
\ref{due}} we have $N^1(X) \cong \ZZ[\xi_{|X}] \oplus \ZZ[F_{|X}]$. For the embedding $X \subset \PP^{N+1}$, the
general hyperplane section (in the countable Zariski topology) $S \subset \PP^N$ of $X$ is an extendable elliptic
surface over $B$ and, as in Example {\rm \ref{uno}}, we have $N^1(S) \cong \ZZ[C] \oplus \ZZ[f]$ by
\cite[Thm.7.5]{moi}. If $L = \O_X(1)$ we get $L_{|S} \eqv C$, $C.f = 4$ and $L_{|S}.f = 4$. With the notation of
Theorem \ref{nonextell}, this is the case $(a, C.f) = (1, 4)$.
\end{exa}

\begin{exa} 
\label{cinque}
Let $\GG(1, 4) \subset \PP^9 = \PP H^0(\GG(1, 4), H)$ be the Grassmannian in its Pl\"ucker embedding $H$ and let $A$ be a very
ample line bundle of degree $v$ on $B$ giving an embedding $B \subset \PP^r = \PP H^0(A)$. Let $N = 10r - 6$ and consider the Segre
embedding $Y = S(B \times \GG(1, 4)) \subset \PP^{N + 5}$. Let $M \cong \PP^{N + 1}$ be a general linear space and let $X = Y \cap M
\subset \PP^{N + 1}$, together with its two projections $p_1 : X \to B$ and $p_2 : X \to \GG(1, 4)$. As in Example {\rm \ref{due}} we
have $N^1(X) \cong \ZZ[p_2^{\ast}H] \oplus \ZZ[F]$, where $F$ is a fiber of $p_1$. For the embedding $X \subset \PP^{N + 1}$, the
general hyperplane section (in the countable Zariski topology) $S \subset \PP^N$ of $X$ is an extendable elliptic surface over $B$
and, as in Example {\rm \ref{uno}}, we have $N^1(S) \cong \ZZ[C] \oplus \ZZ[f]$ by \cite[Thm.7.5]{moi}, where $C =
(p_2^{\ast}H)_{|S}$. Here $L = \O_X(1) \cong p_1^{\ast} A \otimes p_2^{\ast}H$, whence $L_{|S} \eqv C + vf$. Moreover $C.f = L_{|S}.f
= \deg F = \deg \GG(1, 4) = 5$. With the notation of Theorem \ref{nonextell}, this is the case $(a, C.f) = (1, 5)$.
\end{exa}

\begin{exa} \cite[Exa.8.3.9]{ma} {\rm (the example is in fact due to Mori).}
\label{sei}
Let $E$ be a smooth elliptic curve together with a translation $\tau: E \to E$ of order $6$. Let $T$ be the blow-up of
$\PP^2$ at three general points. Then $T$ has an automorphism $\sigma$ of order $6$. Let $X = (T \times E) /
<\sigma, \tau>$ and let $\phi: X \to B = E /<\tau>$ be the natural projection. Then $X$ is a smooth threefold and
$\phi$ is a contraction of an extremal ray (arising from a $(-1)$-curve on $T$) and the fibers $F$ are Del Pezzo
surfaces with $K_F^2 = 6$. Moreover, by \cite[Thm.3.2]{mor}, we have $N^1(X) \cong \ZZ[-K_X] \oplus \ZZ[F]$. We will
prove below that $L = - K_X + h F$ is very ample for $h >> 0$ and that $h^2(\O_X) = 0$. Now a general (in the
countable Zariski topology) $S \in |L|$ is an extendable elliptic surface over $B$ and, as in Example {\rm
\ref{uno}}, we have $N^1(S) \cong \ZZ[C] \oplus \ZZ[f]$ by \cite[Thm.7.5]{moi}, where $C = (-K_X)_{|S}$. Then $L_{|S}
\eqv C + h f$ and $C.f = L_{|S}.f = L^2.F = L_{|F}^2 = K_F^2 = 6$. With the notation of Theorem \ref{nonextell},
this is the case $(a, C.f) = (1, 6)$.
{\rm To see the claim first observe that, since the fibers of $\phi$ are Del Pezzo surfaces, we get, by the Leray
spectral sequence, that $h^2(\O_X) = 0$. Now let $\L$ be any very ample line bundle on $X$, so that there exist
integers $\alpha, \beta$ such that $\L \sim \alpha (-K_X) + \beta F$. Then $\L_{|F} \sim - \alpha K_F$ is ample, so
that $\alpha > 0$. If $j := \lceil \frac{\beta}{\alpha} \rceil$, we get that $- 2 K_X + 2j F \eqv 
\frac{2}{\alpha} \L + 2 (j - \frac{\beta}{\alpha}) F$ is ample, and therefore $H^1(- K_X + h F) = 0$ for $h \geq 2j$
by Kodaira vanishing, since $- K_X + h F = K_X - 2 K_X + h F$. Also $(- K_X + h F)_{|F} = - K_F$ is very ample on
$F$, whence $- K_X + h F$ is base-point free for $h \geq 2j + 1$. We will now prove that $L = - K_X + h F$ is
very ample for $h \geq 2j + 2$. Let $x, y \in X$ be two distinct points. If $x$ and $y$ belong to the same fiber
$F$, we can separate them with sections in $|L|$ since $|L_{|F}|$ is very ample and $H^1(L - F) = 0$. If $x$ and $y$ belong to
two different fibers $F_x$ and $F_y$ respectively, then to separate them just use the fact that $|L - F_x|$ is base-point
free. On the other hand suppose that $x \in X$, $y \in T_x X$ and $d \varphi_L (y) = 0$, where $d \varphi_L$ is the 
differential of the morphism $\varphi_L : X \to \PP H^0(L)$. Arguing as above, $y$ must be tangent to $F_x$, 
contradicting the fact that $|L_{|F_x}|$ is very ample.}
\end{exa}

\begin{exa} 
\label{sette}
Suppose $\E$ has rank $3$, let $Y = \PP \E$ and, for $d = 1, 2$, let $B_i \subset Y, 1 \leq i \leq d$ be sections
of $\pi : Y \to B$ of type $B_i = H_i \cap H'_i$, for general hyperplanes $H_i, H'_i \in |\xi|$. Let
$\varepsilon : X \to Y$ be the blow-up of $Y$ along $B_1, \ldots, B_d$ and denote by $E_1, \ldots, E_d$ the corresponding
exceptional divisors and by $G$ a fiber of $p : X \to B$. As we will see below, the line bundle $L = - K_X + h G$ is very
ample for $h >> 0$ and $G$ is embedded by $L$ as a smooth Del Pezzo surface of degree $9 - d$. Let $S \in |L|$ be a general
hyperplane section (in the countable Zariski topology). We will show that $p_{|S} : S \to B$ is a Weierstrass fibration and
that $\rho(S) = 2 + d$. {\rm To see the assertions claimed above, let us assume that $L$ is very ample. Since $K_S \sim h
G_{|S}$, we have $h^2(\O_S) > h^2(\O_X) = h^2(\O_Y) = 0$ by the birational invariance of $h^{2,0}$. By
\cite[Thm.7.5]{moi} we deduce that $N^1(S) \cong N^1(X)$ has rank $2 + d$. For $1 \leq i \leq d$, we have
$E_i.L.G = - {E_i}_{|G}.{K_X}_{|G} = - {E_i}_{|G}.K_G = 1$, whence $E_i \cap S$ is a section of $p_{|S}$. By the choice of 
$B_1, \ldots, B_d$ we have that each fiber $G$ of $p : X \to B$ is just $\PP^2$ blown-up at $d$ distinct points and $L_{|G}
\sim - K_G$ is very ample, so that $G$ is embedded in $\PP^N = \PP H^0(L)$ by $- K_G$. In particular
$G$ is not ruled by lines and therefore by Castelnuovo-Kronecker's theorem (see for example \cite[LemmaII.2.4]{lo}) we find that $G$
does not have a $(N - 1)$-dimensional family of reducible hyperplane sections. Therefore a general hyperplane $H \in (\PP^N)^{\ast}$
is such that $S = X \cap H$ is smooth and $G \cap H$ is irreducible for all fibers $G$ of $p$. Hence $p_{|S} : S \to B$ does not have
reducible fibers and it must then be a Weierstrass fibration. Finally to see that $L$ is very ample observe that $L_{|G} \sim - K_G$
is very ample by \cite[Thm.V.4.6]{ht}, whence, arguing as in example \ref{sei}, it is enough to find some $h_1 > 0$ such that $H^1(-
K_X + h G) = 0$ for $h \geq h_1$. Then $L = - K_X + h G$ will be very ample for $h \geq h_1 + 2$. To find such $h_1$, for each $b \in
B$ and each fiber $G_b = p^{-1}(b)$, let $G_{b, n} = X \times_B \Spec \O_{B, b} / m_b^n$ be the n-th thickening of $G_b$. Now if $\I$
is the ideal sheaf of $G_b$ we have that $\I^n  / \I^{n+1} \cong \Sym^n \I  / \I^2 \cong \O_{G_b}$, whence, as in \cite[Proof of
Prop.V.3.4]{ht}, there is an exact sequence
\begin{equation}
\label{formal}
0 \hpil \O_{G_b} \hpil \O_{G_{b, n+1}} \hpil \O_{G_{b, n}} \hpil 0.
\end{equation}
Since $H^1(L_{|G_b}) = H^1(- K_{G_b}) = 0$ it follows from \eqref{formal} by induction on $n$ that $H^1(L_{|G_{b, n}}) = 0$ for
each $n \geq 1$ and now the theorem on formal functions \cite[Thm.III.11.1]{ht} gives that $R^1 p_{\ast} L = 0$. Therefore, by the
Leray spectral sequence, we have $H^1(X, L) \cong H^1(B, p_{\ast} L) \cong H^1(B, \pi_{\ast} (\varepsilon_{\ast} L))$. Hence
it will be enough to prove that $H^1(Y, \varepsilon_{\ast} L) = 0$. Now if $g$ is the genus of $B$ and $e$ is the degree
of $\E$ we have that $\varepsilon_{\ast} L \cong \I_{\{B_1 \cup \ldots \cup B_d \}/Y}(3 \xi + \pi^{\ast}(N - K_B - \det \E))$ for
some line bundle $N \in \Pic^h B$. Now let
$M = 3 \xi + \pi^{\ast}(N - K_B - \det \E) \eqv 3 \xi + (h - 2g + 2 - e) F$ and consider the exact sequence
\begin{equation}
\label{ideal}
0 \hpil \I_{\{B_1 \cup \ldots \cup B_d \}/ Y} (M) \hpil M \hpil \O_{B_1 \cup \ldots \cup B_d} (M) \hpil 0. 
\end{equation} 
It is easily seen that, for $h >> 0$ we have $H^1(M) = 0$, whence, from \eqref{ideal}, it remains to show that $H^0(M) \to
H^0(\O_{B_1 \cup \ldots \cup B_d} (M))$ is surjective. For $d = 1$ the required surjectivity follows easily by the definition of
$B_1$. For $d = 2$ we need to show, for $\{i, j \} = \{1, 2 \}$, that the maps $H^0(\I_{B_i/Y}(M) \to H^0(\O_{B_j}(M))$ are
surjective. To this end consider the following exact diagram, defining the sheaf $\F$, where the middle exact sequence is the Koszul
resolution of $B_j \subset Y$:
\begin{equation*}  
\xymatrix{& 0 \ar[d] & 0 \ar[d] & 0 \ar[d] & & \\
0 \ar[r] & \I_{B_i / Y}(M - 2\xi)  \ar[r] \ar[d] & \I_{B_i / Y}(M - \xi)^{\oplus 2} \ar[r] \ar[d] & \F \ar[r] \ar[d] & 0 \\
0 \ar[r] & M - 2 \xi \ar[r] \ar[d] & (M - \xi)^{\oplus 2} \ar[r] \ar[d] & \I_{B_j / Y}(M) \ar[r] \ar[d] & 0
\\ 0 \ar[r] & (M - 2 \xi)_{|B_i} \ar[r] \ar[d] & (M - \xi)_{|B_i}^{\oplus 2} \ar[r] \ar[d] & M_{|B_i} \ar[r] \ar[d] & 0 \\ 
& 0 & 0 & 0 & & }
\end{equation*} 
Now we just need $H^1(\F) = 0$, which, in turn, follows from $H^1(\I_{B_i / Y}(M - \xi)) = H^2(\I_{B_i / Y}(M - 2 \xi)) = 0$ for
$i = 1, 2$. Finally the latter two vanishings follow easily from the Koszul resolution of $B_i \subset Y$.} 
\end{exa}

We end the section with two simple examples of extendable elliptic $K3$ surfaces of rank two.

\begin{exa} 
Consider $S_4 \subset \PP^3$ a general quartic surface containing a line $C$. The pencil of planes through $C$
gives a fibration $\pi : S \to \PP^1$ by elliptic plane cubic curves, $\rho(S) = 2$ (by \cite[Cor.II.3.8]{lo} or by
\cite[Thm.1.1]{kn}), $C^2 = -2$, $C.f = 3$ and $S_4 \subset \PP^3$ is extendable.
\end{exa}

\begin{exa}
Consider $S_{2,3} \subset \PP^4$ a general complete intersection containing a linearly normal elliptic quintic $C
\subset \PP^4$. Then $\rho(S) = 2$ (by \cite[Thm.III.2.1]{lo} or by \cite[Thm.1.1]{kn}), and $S$ has an
elliptic fibration $\pi : S \to \PP^1$ given by $f \sim 5H - 3C$, $C^2 = 0$, $C.f = 25$ and $S_{2,3} \subset \PP^4$ is
extendable.
\end{exa}

As a matter of fact, if $\pi : S \to \PP^1$ is a $K3$ elliptic surface with $\rho(S) = 2$ then, by
\cite[Thm.5.1]{ko} and \cite[Exa.1.4.33]{la2}, there exists an irreducible effective curve $C \subset S$ such that
$\overline{NE}(S) = \{\alpha C + \beta f \mid \alpha, \beta \geq 0\}$ such that either $C^2 = -2$ and $C
\cong \PP^1$ or $C^2 = 0$ and $p_a(C) = 1$.

\section{A non extendability condition for other fibered surfaces}
\label{other}

Following the ideas of Section \ref{eoes} we present a criterion for a fibered surface $S\subset\PP^N$ not to be
extendable.

\begin{prop}
\label{criterio}
Let $Y \subset \PP^N$ be a smooth irreducible nondegenerate variety. Let $p : Y \to Z$ be a surjective morphism
onto a projective variety $Z$ having a finite-to-one morphism into an abelian variety and suppose that $\dim Y
\geq \dim Z + 1$. If the general fiber $f \subset \PP^N$ of $p$ is not extendable and is not a linear subspace,
then $Y$ is not extendable to a normal variety $X \subset \PP^{N + 1}$ with rational singularities.
\end{prop}

\begin{proof}
Let $X \subset \PP^{N+1}$ be a normal variety with rational singularities containing $Y \subset \PP^N$ as
a hyperplane section. By \cite[Thm.5.2.3]{bs} the morphism $p: Y \to Z$ extends to a morphism $\overline{p} : X \to
Z$. Let $f$ and $F$ denote their respective fibers. Since $f \subset \PP^N$ is not extendable we have that $F
\subset \PP^{N + 1}$ is a cone over $f$, whence $F$ is singular since $f$ is not a linear subspace. 

On the other hand $F$ is smooth by Bertini's theorem and this contradiction proves the theorem.
\end{proof}

We have the following nice consequence of Proposition \ref{criterio} (by Proposition \ref{1.10} and Zak's theorem \cite[page
277]{z}). 

\begin{cor}
Let $S$ be a smooth irreducible surface with a surjective morphism $\pi : S \to B$ onto a smooth
irreducible curve with $g(B) > 0$. Let $L$ be a very ample line bundle on $S$ such that, on a general fiber
$f$ of $\pi$ we have that $H^0(S, L) \to H^0(f, L_{|f})$ is surjective, $g (f) > 0$ and the Gaussian map $\Phi_{L_{|f},
\omega_f}$ is surjective. Then, in the embedding $S \subset \PP H^0(L)$, $S$ is not extendable to a normal variety
with rational singularities.
\end{cor}

For example one can take any fibration whose general fiber is not trigonal and not isomorphic to a plane quintic
and line bundles $L = 2K_S + f + A$ for any line bundle $A$ on $S$ such that $K_S + A$ is big and nef and $A.f
\geq 1$. In this case $H^1(L - f) = 0$ by Kawamata-Viehweg vanishing and the Gaussian map $\Phi_{L_{|f},
\omega_f}$ is surjective by Theorem \ref{Thm.2} (since $\deg L_{|f} = 2K_S.f + A.f \geq 4g(f) - 3$).

\end{document}